\documentclass[preprint,authoryear,sort&compress,11pt]{elsarticle}
\pdfoutput=1
\usepackage{fullpage}
\usepackage[utf8]{inputenc}
\usepackage[english]{babel}
\usepackage{graphicx,amsmath,amssymb, amsthm, multicol, bm}
\usepackage[table]{xcolor}
\usepackage[section]{placeins}
\usepackage{tikz-cd}[table,xcdraw]
\usepackage{colortbl}
\usetikzlibrary{shapes.geometric}
\usepackage[shortlabels]{enumitem}
\usepackage{soul}
\usepackage{txfonts}
\usepackage{tocbibind}
\usepackage{pifont}
\usepackage{float}
\usepackage{calc}
\usepackage{tikz} 
\usepackage{adjustbox}
\usetikzlibrary{decorations}
\usepackage[hidelinks,unicode,verbose,hypertexnames=false]{hyperref}
\usepackage[normalem]{ulem}

\journal{European Journal of Operational Research}
\makeatletter
\def\ps@pprintTitle{%
     \let\@oddhead\@empty
		 \let\@evenhead\@empty
     \def\@oddfoot
       {\hbox to \textwidth%
        {\ifnopreprintline\relax\else
        \@myfooterfont%
         \ifx\@elsarticlemyfooteralign\@elsarticlemyfooteraligncenter%
           \hfil\@elsarticlemyfooter\hfil%
         \else%
         \ifx\@elsarticlemyfooteralign\@elsarticlemyfooteralignleft%
           \@elsarticlemyfooter\hfill{}%
         \else%
         \ifx\@elsarticlemyfooteralign\@elsarticlemyfooteralignright%
           {}\hfill\@elsarticlemyfooter%
         \else%
            \begin{minipage}[t]{\textwidth}   Preprint accepted in \ifx\@journal\@empty%
                 Elsevier%
            \else\@journal\fi\hfill August 10, 2024\\
						Published version available at \rurl{10.1016/j.ejor.2024.08.009}\medskip\\ 
						\begin{minipage}[b]{0.7\textwidth}\textcopyright\ 2024. This work is licensed under a Creative Commons\\
			      Attribution-NonCommercial-NoDerivatives
						4.0 International License\end{minipage} 
							\hfill\includegraphics[width=0.2\textwidth]{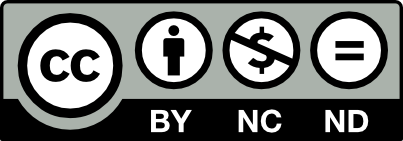}
						\end{minipage}\fi%
         \fi%
         \fi%
         \fi%
         }
       }%
     \let\@evenfoot\@oddfoot}
\makeatother
\newcommand\rurl[1]{%
  \href{http://doi.org/#1}{\nolinkurl{#1}}%
}

\newlength{\tailarrowtiplength}
\newlength{\tailarrowtipwidth}
\newlength{\tailarrowtailwidth}
\newlength{\tailarrowjointwidth}

\pgfdeclaredecoration{tail arrow decoration}{initial}{
    \state{initial}[width=\pgfdecoratedpathlength, next state=final] {
        \pgfpathlineto{\pgfpoint{0pt}{0.5\tailarrowtailwidth}}
        \pgfpathlineto{\pgfpointadd{\pgfpointdecoratedpathlast}{\pgfpoint{-1\tailarrowtiplength}{0.5\tailarrowjointwidth}}}
        \pgfpathlineto{\pgfpointadd{\pgfpointdecoratedpathlast}{\pgfpoint{-1\tailarrowtiplength}{0.5\tailarrowtipwidth}}}
        \pgfpathlineto{\pgfpointdecoratedpathlast}
        \pgfpathlineto{\pgfpointadd{\pgfpointdecoratedpathlast}{\pgfpoint{-1\tailarrowtiplength}{-0.5\tailarrowtipwidth}}}
        \pgfpathlineto{\pgfpointadd{\pgfpointdecoratedpathlast}{\pgfpoint{-1\tailarrowtiplength}{-0.5\tailarrowjointwidth}}}
        \pgfpathlineto{\pgfpoint{0pt}{-0.5\tailarrowtailwidth}}
        \pgfpathclose
    }
    \state{final} {
        \pgfpathmoveto{\pgfpointdecoratedpathlast}
    }
}

\tikzset{
    tail arrow tip length/.code={
        \setlength{\tailarrowtiplength}{#1}
    },
    tail arrow tip width/.code={
        \setlength{\tailarrowtipwidth}{#1}
    },
    tail arrow tail width/.code={
        \setlength{\tailarrowtailwidth}{#1}
    },
    tail arrow joint width/.code={
        \setlength{\tailarrowjointwidth}{#1}
    },
    tail arrow tip length/.default={ 2em },
    tail arrow tip width/.default={ 2em },
    tail arrow tail width/.default={ 1em },
    tail arrow joint width/.default={ .2em },
    tail arrow/.style={
        tail arrow tip length,
        tail arrow tip width,
        tail arrow tail width,
        tail arrow joint width,
        decorate,
        decoration={tail arrow decoration}
    },
    vertup/.style={anchor=south, rotate=90, inner sep=1.5mm},
    vertdn/.style={anchor=south, rotate=90, inner sep=1.5mm}
}

\newcommand*{\mask}[2]{%
    \makebox[\widthof{\(#2\)}]{\(#1\)}%
}

\newtheorem{theorem} {Theorem}[section]

\newtheorem{lemma}[theorem] {Lemma}

\newtheorem{definition}[theorem] {Definition}
\newtheorem{corollary}[theorem]{Corollary}
\theoremstyle{remark}
\newtheorem{example}[theorem] {Example}
\newtheorem{pozn}[theorem]{Remark}

\newcommand{\DMU}{\mathrm{DMU}}
\newcommand{\ones}{\bm{e}}
\newcommand{\T}{\mathcal{T}}
\newcommand{\F}{\mathcal{F}}
\newcommand{\HP}{\mathcal{H}}
\newcommand{\dom}{\operatorname{dom}}
\newcommand{\im}{\operatorname{im}}
\newcommand{\cl}{\operatorname{cl}}
\renewcommand{\int}{\operatorname{int}}

\newcommand{\R}{\mathbb{R}}
\newcommand{\m}{\mathrm{min}}
\newcommand{\M}{\mathrm{max}}
\newcommand{\X}{\bm{X}}
\newcommand{\Y}{\bm{Y}}
\newcommand{\x}{\bm{x}}
\newcommand{\y}{\bm{y}}
\newcommand{\z}{\bm{z}}

\newcommand{\g}{\bm{g}}
\newcommand{\p}{\bm{\phi}}
\newcommand{\la}{\bm{\lambda}}
\newcommand{\0}{\bm{0}}
\newcommand{\s}{\bm{s}}
\newcommand{\ab}{\bm{a}}
\newcommand{\bb}{\bm{b}}
\newcommand{\e}{\bm{e}}
\newcommand{\D}{\mathcal{D}}
\newcommand{\uu}{\bm{u}}
\newcommand{\vv}{\bm{v}}

\makeatletter
\newcommand{\myitem}[1]{%
\item[#1]\protected@edef\@currentlabel{#1}%
}
\makeatother
\makeatletter
\newcommand\newtag[2]{#1\def\@currentlabel{#1}\label{#2}}
\makeatother

\renewenvironment{thebibliography}[1]{%
\begin{oldthebibliography}{#1}%
\setlength{\baselineskip}{1.2em}
\linespread{1.2}
\setlength{\parskip}{.3ex}%
\setlength{\itemsep}{.3em}%
}%
{%
\end{oldthebibliography}%
}

\begin{document}

\begin{frontmatter}

\title{
On indication, strict monotonicity, and efficiency of projections\\ in a general class of path-based data envelopment analysis models}
\author[fmfi]{Margar\'eta Halick\'a}
\ead{halicka@fmph.uniba.sk}
\author[fmfi]{M\'aria Trnovsk\'a\texorpdfstring{\corref{cor1}}{}}
\ead{trnovska@fmph.uniba.sk}
\author[bayes]{Ale\v{s} \v{C}ern\'{y}}
\ead{ales.cerny.1@city.ac.uk}

\cortext[cor1]{Corresponding author}
\address[fmfi]{Faculty of Mathematics,  Physics and Informatics, Comenius
 University in Bratislava,  Mlynsk\'a dolina, 842 48 Bratislava, Slovakia}
 \address[bayes]{Bayes Business School, City St George's, University of London, 106 Bunhill Row, London EC1Y 8TZ, UK}

\begin{abstract} 
Data envelopment analysis (DEA) theory formulates a number of desirable properties that DEA models should satisfy. Among these,  indication, strict monotonicity, and strong efficiency of projections tend to be grouped together in the sense that, in individual models, typically, either all three are satisfied or all three fail at the same time. Specifically, in slacks-based graph models, the three properties are always met; in path-based models, such as radial models,  directional distance function models, and the hyperbolic function model, the three properties, with some minor exceptions, typically all fail.

Motivated by this observation, the article examines relationships among indication, strict monotonicity, and strong efficiency of projections in the class of path-based  models over variable returns-to-scale technology sets. Under mild assumptions, it is shown that the property of strict monotonicity and strong efficiency of projections are equivalent, and that both properties imply indication. This paper also characterises a narrow class of technology sets and path directions for which the three properties hold in path-based models.
\end{abstract}

\end{frontmatter}

\section{Introduction}

Data envelopment analysis (DEA) is a non-parametric method measuring relative efficiency within a group of homogeneous units that use multiple inputs to produce multiple outputs.
This paper focuses on DEA models that integrate an analytical description of the technology set and an efficiency measure into a single mathematical optimisation programme. 
The output of the programme for the unit under evaluation is its efficiency score and a certain benchmark/projection (a point on the frontier of the technology set) from which the score is derived.  The concepts of \textit{efficiency score} and \textit{projection} enter the formulations of the so-called desirable properties of the models.
    
The most comprehensive lists of desirable properties
are provided in \cite*{cooper.al.99} and \cite{sueyoshi.sekitani.09.ejor.196}. These properties include indication of strong efficiency; homogeneity;  strict/weak monotonicity; boundedness; unit invariance; and translation invariance. The selected DEA models are then classified on the basis of these criteria. In the works \cite{halicka.trnovska.21}
and \cite*{halicka2024unified}, the property of strongly efficient projection is added to the desirable properties and the fulfilment of all properties is studied within two wide classes of DEA models. It is observed that the three properties: indication, strong efficiency of projection, and strict monotonicity usually appear as a trio, either all three are fulfilled or not fulfilled in the given model.
  
This article examines the connection between the above mentioned three properties, the meaning and importance of which are as follows. \\
 {\bf Indication:} 
 \textit{The efficiency score is equal to one if and only if the evaluated unit is strongly efficient.} This property allows revealing the strong efficiency of a unit just from the mere value of its efficiency score.  
 \\
 {\bf Strict monotonicity:} \textit{An increase in any input or decrease of any output relative to the evaluated unit, holding other inputs as well as outputs constant, reduces the efficiency score.} This property states that the measure of efficiency is fair -- if a unit dominates another unit, then the former one has greater efficiency score. In the work \cite{pastor.al.99.enhanced}, this property is interpreted as the sensitivity of the measure to changes in inputs or outputs.
\\
 {\bf Strong efficiency of projections:} \textit{ Projections generated by the model are strongly efficient.} This property is sometimes called \emph{the property of efficient comparison}. It states that the score is based on the comparison of the evaluated unit with a strongly efficient one, and therefore the value of the efficiency score is not overestimated. It could be viewed as an extension of the indication property (dealing with units with efficiency score equal to one) to any unit, and it could be alternatively formulated as: the efficiency score accounts for all sources of inefficiency if and only if the projection point is strongly efficient.  

Note that the desirable properties of the measures over general technology sets were already formulated and analysed prior to the emergence of DEA -- in the framework of economic production theory.\footnote{In the economic production theory, the measures of technical efficiency or inefficiency are called indices and the desirable properties of the indices are called axioms.}
The first work that formulated certain properties that an input (or output)-based efficiency measure should satisfy was the work of \cite{fare.lovell.78}. In this work, four properties were formulated, which, in addition to the three mentioned above, also included the property of homogeneity. 
Next, in \cite{russell1985measures} it was shown that, for a special type of measure, the efficient comparison property is redundant with respect to the other three properties. Certain objections were also raised regarding the unclear formulation of this property, worded as: \textit{``\textellipsis comparison to efficient input vectors (the measure compares each feasible input vector to an efficient input vector).''}

Later, in both areas, the list of desirable properties was expanded and the properties were examined not only in connection with input or output models, but also for graph measures of efficiency such as hyperbolic or additive measures (e.g. \citealp{pastor.al.99.enhanced} and \citealp{russell.schworm.11}). In these and other works, the efficient comparison property was no longer mentioned, objections to the vagueness of the definition of this property were not reconsidered, and the validity of the conclusions of \cite{russell1985measures} about its redundancy for other types of measures was not verified. 

In actual fact, DEA does provide mathematical tools for defining the projection and also supplies additional reasons to not neglecting the property of efficient comparison.
Thanks to DEA, the models can be formulated as mathematical programming problems, and the relevant benchmarks/projections can be easily identified from their optimal solutions. This makes the efficient comparison property well-defined. Furthermore,  the efficient comparison property (as well as the indication property) can be easily verified in the DEA models for each of the assessed units using the standard second-phase method. 
Nevertheless, the efficient comparison property is only seldom included in the list of desirable properties,%
\footnote{To the best of our knowledge, the efficient comparison property is explicitly listed among desirable properties of DEA models only in \cite{halicka.trnovska.21}, \cite{halicka2024unified} and \citet[p.~40]{pastor2022benchmarking}, where in the last cited work this property appears as an extension of the standard `indication property' under (E1b).}
and instead many authors settle for the weaker property of indication. 
However, the efficient comparison property and its verification in specific models are particularly important in DEA. Namely, DEA approximates the most widely used variable returns-to-scale technology using observed data coupled with the postulates of convexity and free disposability. This leads to polyhedral technologies whose large parts of the frontiers are not strongly efficient. Then, individual models projecting units on the frontier may derive scores from projections that are weakly but not strongly efficient.

There are two basic ways of searching for the benchmarks in DEA, leading to the classification of models into path-based models and slacks-based models introduced in the work \cite{russell.schworm.18}.
As explained there, the slacks-based measures are expressed in terms of additive or multiplicative slacks for all inputs and outputs, and particular measures are generated by specifying the form of aggregation over the coordinate-wise slacks. On the other hand, the path-based
measures are expressed in terms of a common contraction/expansion factor, and particular measures are generated by
specifying the parametric path leading from the assessed unit towards the frontier of the technology set. In this class, the projection is uniquely determined as the point at which the path leaves the technology set. However, the slacks-based models may provide multiple benchmarks and hence it may not be apparent from which benchmark the efficient score is derived.

The desirable properties of slacks-based models over variable returns-to-scale (VRS) technology sets were analysed in \cite{halicka.trnovska.21} through a general scheme. The scheme encompassed all commonly used models, such as the Slacks-Based Measure (SBM) model of \cite{tone2001slacks}, the Russell Graph Measure model of \cite*{fare.al.85}, the Additive Model (AM) of \cite{charnes.al.85}, and the Weighted Additive Models (WAM) including the Range Adjusted Measure (RAM) model  (\citealp*{cooper.al.99}) and  the Bounded Adjusted Measure (BAM) model (\citealp{cooper2011bam}). It was shown in \cite{halicka.trnovska.21}  that the scores of these models are derived from benchmarks that may not be unique but are all strongly efficient, and the resulting score does not depend on the choice of the benchmark. As a consequence, all models in this class satisfy the indication and efficient comparison properties and, therefore, they all account for all sources of inefficiencies. The strict monotonicity property had to be individually assessed for models belonging to this scheme and was proven for all standard models with the exception of the BAM model, which only met the weak monotonicity.

Quite different results were obtained in the class of path-based models under VRS technology sets. These models include the radial BCC input or output-oriented models (\citealp*{banker.al.84}),  the directional distance function (DDF) model (\citealp*{chambers1996benefit,chambers.al.98}), and the hyperbolic distance function (HDF) model (\citealp*{fare.al.85}). 
\cite{halicka2024unified} analysed this class through a general scheme that depends on the parameters choices (convex functions and directional vectors) and covers all standard path-based models. 
 The scheme allows for negative data. 

The results showed that most models fail simultaneously all three properties of indication, strict monotonicity, and efficient comparison. However, the authors also presented very specific examples of path-based models that, in the case of technology sets with one input and one output, met all three properties.
The simultaneous failure or success of the three properties invites the following questions,
\begin{enumerate}[(i)]
    \item to what extent the three properties are related;
    \item can one find sufficient conditions under which path-based models satisfy all three properties.
\end{enumerate}
The aim of this article is to provide answers to these two questions.

With regard to the first question, it is fairly immediate that the property of efficient comparison implies the indication property; such an implication is also valid outside of the scheme of path-based models. However, the reverse implication is not apparent, and we will show that it is not even generally valid.  Moreover, the connection between the efficient comparison property and the property of strict monotonicity is also unclear. In the framework of the general scheme of path-based models, we will show the equivalence of these two properties under mild assumptions. To the best of our knowledge, this equivalence will be identified and proved here for the first time. 

Obtaining an answer to the second question can be important for practitioners --- is it possible to attain the three desirable properties by a suitable choice of parameters of path-based models?  An analysis outlined in \cite{halicka2024unified}  indicates that appropriately modified range directions of \cite{portela.al.04} tend to improve the properties of the model.  These directions work well in single-input and single-output examples but not in all three-dimensional examples. More extensive numerical experiments with two inputs and two outputs carried out in \cite{halicka2024unified} show that these directions can significantly reduce the number of units not projected onto the strongly efficient frontier but cannot eliminate the presence of such units completely in general. Therefore, it seems that the problem partly stems from the characteristics of the VRS technology set generated by the data.
In the second part of this article, we confirm this conjecture. We provide  a characterisation of type of data, or equivalently, the type of technological set that, together with the aforementioned choice of directions, will ensure the three properties are met. We term this type of technology the \emph{ideal technology}.  We also offer a practical recipe to determine whether given data generates an ideal technology set.

The paper is organised as follows. The preliminaries required to accurately address the topics discussed in this article are outlined in Section~\ref{S2}. Section~\ref{S3} examines the connections among the three desirable properties. In Section~\ref{S4}, particular directions and technology sets are outlined to ensure that all three properties are satisfied. Section~\ref{S5} contains several numerical illustrations on real data. The concluding Section~\ref{S6} summarises our theoretical findings and draws lessons for practitioners. \ref{A} analyses the facial structure of a general VRS technology set, \ref{B} gives equivalent characterizations of the ideal technology, and \ref{C} contains the proof of Lemma~\ref{AcorPW}. 
\section{Preliminaries}\label{S2}

Let  $\R^d$ be  the $d$-dimensional Euclidean space and $\R^d_+$ its non-negative orthant. Bold lowercase letters denote column vectors, and bold uppercase letters denote matrices.
The superscript $^\top$ denotes the transpose of a column vector or a matrix.
For a vector $\z_o\in \R^d$, $z_{ko}$ denotes its $k$-th component, and hence $\z_o=[z_{1o},\dots,z_{do}]^\top$.  The symbol $\ones$ denotes a vector of ones. 

We will consider a production process with $m$ inputs and $s$ outputs. For any two input-output vectors $(\x_o,\y_o),(\x_p,\y_p)\in\R^m\times\R^s$  we will use the notation
\begin{itemize}
    \item $(\x_o,\y_o)\succsim (\x_p,\y_p)$ if $(\x_o, \y_o)$ \emph{weakly dominates}  $( \x_p, \y_p)$, i.e., $ \x_o\le \x_p$ and $ \y_o\ge\y_p$;
    \item $(\x_o,\y_o)\succnsim (\x_p,\y_p)$ if $(\x_o, \y_o)$ \emph{dominates}  $( \x_p, \y_p)$, i.e,   $ \x_o\le \x_p$, $ \y_o\ge\y_p$, and $(\x_o,\y_o)\neq (\x_p,\y_p)$;
    \item $(\x_o,\y_o)\succ (\x_p,\y_p)$ if $( \x_o, \y_o)$ \emph{strictly dominates}  $( \x_p, \y_p)$, i.e.,    $ \x_o<\x_p$ and $ \y_o> \y_p$.
\end{itemize}

\subsection{Technology set}\label{SS2.1}

Consider a set of $n$ decision-making units ${\DMU}_j$ ($j=1,\dots, n$) with observed input--output vectors $(\x_j,\y_j)\in \R^m\times\R^s$. 
The input--output data of $(\DMU_j)_{j=1}^n$ are arranged into  the $m\times n$ input and $s\times n$ output matrices $\X=[\x_1,\dots,\x_n]$  and $\Y=[\y_1,\dots ,\y_n]$,  respectively.
No assumption about the non-negativity of the data is made at this point. The non-negativity requirement may follow later from other assumptions placed on the models.

Based on the given data  $\X, \Y$ we consider the  technology set
\begin{equation}\label{T}
\T=\left\{(\x,\y)\in\R^m\times \R^s\ | \
\X\la\le \x,\ \Y\la\ge \y,\ \la\ge \0, \ \ones^\top\la =1
\right\},
\end{equation}
corresponding to variable returns to scale (VRS).  Note that
the common non-negativity of $(\x,\y)$ is not imposed here. It follows from \eqref{T} that the set $\T$ is closed, has a non-empty interior (denoted $\int\T$) and its boundary $\partial\T$ satisfies $\partial \T=\T\setminus \text{int}\T$. Elements of ${\T}$ will be called units. By $(\x_o,\y_o)$ we denote a unit from ${\T}$ to be evaluated. 

The point $(\x^{\m},\y^{\M})\in \R^m\times\R^s$, with elements $x^{\m}_i=\min_jx_{ij}$, for $i=1,\dots m$ and $y^{\M}_r=\max_jx_{rj}$ for $r=1,\dots s$, respectively, is called the  \emph{ideal point of $\T$} in DEA (see, e.g., \citealp{portela.al.04}). The ideal point typically does not belong to $\T$, in which case it dominates every unit in $\T$. The technology set, for which $(\x^{\m},\y^{\M})\in \T$, will be called \emph{trivial}. Clearly, a trivial technology is an affine transformation of a non-negative orthant. 

 A unit $( \x_o, \y_o)\in \T$ is called \emph{strongly efficient}%
  \footnote{This is the well known Pareto--Koopmans efficiency.   Some authors call such units Pareto efficient, or fully efficient - see the discussion in \citet[p. 45]{cooper.al.07}.}
if no other unit in $\T$ dominates  $( \x_o, \y_o)$, i.e.,  if  the property that $( \x, \y)\in \T$  dominates $( \x_o, \y_o)$ yields  $( \x, \y)=( \x_o, \y_o)$. A unit $(\x_o,\y_o)\in \T$ is called \emph{weakly efficient} if there is no unit in $\T$ that strictly dominates $(\x_o,\y_o)$.
  
Evidently, any strongly efficient unit is weakly efficient, and the weakly efficient units lie on the boundary of the technology set. The converse is also true: every unit on the boundary $\partial \T$ is weakly efficient because the definition of $\T$ in \eqref{T} does not impose the non-negativity assumption on the units therein. 
The boundary $\partial \T$ is thus uniquely partitioned into the \emph{strongly efficient frontier} $\partial^S\T$ containing all strongly efficient units and the remaining part $\partial^W\T:= \partial\T\setminus \partial^S\T$, which consists of the weakly but not strongly efficient units. In this paper we will refer to the remaining part of the boundary as the \emph{weakly efficient frontier}. Thus, we have $\partial \T=\partial^S\T\cup\partial^W\T$ and $\partial^S\T\cap\partial^W\T=\emptyset$. 

\subsection{A general scheme for path-based models}\label{SS2.2}

We now recall the general scheme (GS) for path-based models from \cite{halicka2024unified}.  The scheme depends on both the choice of a prescription $\g$ that defines the directional vector $\g_o=(\g_o^x,\g_o^y)\gneqq 0$ for each $(\x_o,\y_o)$,
and the choice of real functions $\psi^x$ and $\psi^y$ that together with their
domains ($\dom$) and images ($\im$) satisfy the following assumptions:
\begin{enumerate}[label={\rm(A\arabic{*})}, ref={\rm(A\arabic{*})}]
 \item\label{A1}
 $\dom(\psi^x)=(a^x,\infty)$ with $a^x\in  \{-\infty,0\}$ and $\dom(\psi^y)=(a^y,\infty)$ with $a^y\in  \{-\infty,0\}$;
 \item\label{A1bis} $\psi^x$ is smooth, concave,  increasing and $\psi^y$ is smooth, convex,  decreasing;
 \item\label{A2}  $\psi^x(1)= \psi^y(1)=1$;
 \item\label{A3}  $\text{im}(\psi^x)=(b^x,\infty)$ with  $b^x=-\infty$ if $\g^y=\0$ and $b^x\in \{-\infty,0\}$ otherwise; $\im(\psi^y)=(b^y,\infty)$ with $b^y\in \{-\infty,0\}$.
 \end{enumerate}
 The scheme is built on the technology set $\T$ specified in \eqref{T} and is consequently also dependent on the matrices $\X,\Y$. 
 Any choices of $(\X,\Y, \g,\psi^x,\psi^y)$ that satisfy the requirements mentioned above will be referred to as \emph{admissible parameters} for a model in the GS scheme.
 
For a fixed choice of admissible parameters $(\X,\Y, \g,\psi^x,\psi^y)$, the (path-based) GS  model 
for assessment of $(\x_o,\y_o) \in \T$ is defined by
\begin{subequations}\label{general}
\begin{align}
(\text{GS})_o\qquad\min\ &{} \theta \label{general1}\\
     &{} \X\la \le \x_o + (\psi^x(\theta)-1) \g_o^x,\label{general2}\\
     &{} \Y\la \ge \y_o + (\psi^y(\theta)-1)\g_o^y,\label{general3}\\
     &{} \ones^\top\la =1, \quad      \la \ge \0.\label{general4}
\end{align}
\end{subequations}

 The right-hand sides of \eqref{general2} and \eqref{general3}, denoted by
\begin{equation}\label{pathdef}
\p_o^x(\theta):=\x_o + (\psi^x(\theta)-1)\g_o^x \quad\text{and}\quad
\p_o^y(\theta):=\y_o + (\psi^y(\theta)-1)\g_o^y,
\end{equation}
 define a smooth path $\p_o(\theta)= (\p_o^x(\theta),\p_o^y(\theta))$ in the input-output space $\R^m\times \R^s$ parametrized by $\theta\in {\cal{D}}=\dom(\psi^x)\cap\dom(\psi^y)$. 

For any choice of admissible parameters and each $(\x_o,\y_o)\in \T$,
the well-definedness of the programme $(GS)_o$ and other useful properties of the path $\p_o$ are established in Theorem 3.1 of \cite{halicka2024unified}. According to this theorem, the path $\p_o$ passes through the point $(\x_o,\y_o)\in \T$ at $\theta=1$ (i.e., $\p_o(1)=(\x_o,\y_o)$), and for decreasing values of $\theta$ it moves towards the boundary of $T$ gradually passing through points that dominate one another. This property of the path can be formally expressed as $\p_o(\theta_1)\succnsim \p_o(\theta_2)$ for $\theta_1<\theta_2$,  and we will refer to it as \emph{the monotonicity of the path with respect to $\theta$}. Finally, the path leaves $\T$ at some $\p_o(\theta^*_o)\in\partial \T$, where  $\theta^*_o\le 1$, and  $\theta^*_o$ is the optimal value in $(GS)_o$.
The optimal value $\theta^*_o$ is called the \emph{efficiency score}, or alternatively, the \emph{value of the efficiency measure} for $(\x_o,\y_o)$. The point $(\p_o^x(\theta_o^*),\p_o^y(\theta_o^*))$ on the path $\p_o$ is called the \emph{projection} of $(\x_o,\y_o)$ in the GS model. 
 
\subsection{Standard path-based models}\label{SS2.3}
It is easy to see that the well-known BCC input and output models \citep{banker.al.84}, the hyperbolic distance function model (HDF)  \citep{fare.al.85} as well as the generic directional distance function model (DDF) \citep{chambers1996benefit,chambers.al.98} can be equivalently rewritten in the form of the general scheme \eqref{general}. The scheme also includes the so-called generalised distance function model (GDF) introduced by \cite{chavas1999generalized}.
These models, taken in conjunction with the usual assumptions on the positiveness of the data, will be called the \emph{standard path-based models}. The corresponding parameterizations are shown in Table~\ref{Tab:values}. 
\begin{table}[ht!]
\centering
 \begin{tabular}{l c c }
 \hline\\[-1.8ex]
 Model   &$\p_o^x(\theta)$ &  $\p_o^y(\theta)$  \\[0.8ex]
 \hline\\[-1.8ex]
 BCC-I \cite{banker.al.84}  & $\x_o+(\theta-1) \x_o$  & $\y_o$ 
  \\[0.2ex]
 BCC-O \cite{banker.al.84}  & $\x_o$ & $\y_o+(\frac{1}{\theta}-1)\y_o$   \\[0.4ex]
 DDF-g \cite{chambers1996benefit,chambers.al.98} &  $\x_o+(\theta-1)\g_o^y$ & $\y_o+(2-\theta-1)\g_o^y$  \\[0.4ex]
 HDF \cite{fare.al.85}& $\x_o+(\theta-1) \x_o$ & $\y_o+(\frac{1}{\theta}-1)\y_o$  \\[0.4ex]
 GDF \cite{chavas1999generalized} $p\in [0,1]$& $\x_o+(\theta^{1-p}-1) \x_o$ & $\y_o+(\theta^{-p}  -1)\y_o$   \\[0.4ex]
 \hline
 \end{tabular}
 \caption{Parameterization of the standard path-based models.}
\label{Tab:values}
\end{table}

 The choice $g^x_o=0$ or $g^y_o=0$ for all $(\x_o,\y_o) \in \T$  leads to output or input oriented models, respectively. Among the standard path-based  models, only DDF is formulated with general directional vectors. The BCC as well as the HDF use $g^x_o=\x_o$ and / or $g^y_o=\y_o$.  A generalisation of HDF towards general directions is introduced in \cite{halicka2024unified}. 
     
\subsection{Properties of the general model}\label{SS2.4}

 In \cite{halicka2024unified}, the GS models are analysed in light of ten desired properties. Three of these merit further investigation. First, we recall
the precise definitions of \cite{halicka2024unified} abridged to suit our needs here.

 \begin{itemize}
\myitem{(ID)}\label{P2a} 
{\bf Identification of strong efficiency.} 
If for a given $(\x_o,\y_o)\in\T$ one has $\theta_o^*=1$, then $(\x_o,\y_o)\in\partial^S\T$. 
\myitem{(PR)}\label{P3} 
{\bf Strong efficiency of projections.}  One has 
$\phi_o(\theta_o^*)\in \partial^S\T$ for each $(\x_o,\y_o)\in\T$. 
\myitem{(MO)}\label{P8} 
{\bf Strict monotonicity.} If  $(\x_o,\y_o)\succnsim (\x_p,\y_p)$ for some $(\x_o,\y_o), (\x_p,\y_p)\in \T$, then $\theta_o^*<\theta^*_p$. 
\end{itemize}

\begin{pozn} Note that (ID) is the \lq only if \rq  part of the property known in DEA as Indication: $\theta_o^*=1$, if and only if $(\x_o,\y_o)\in\partial^S\T$.  
Unlike the \lq only if \rq  part of Indication, the  \lq if \rq part is universally satisfied by each model of the GS scheme (see Theorem 2 in \citealp{halicka2024unified}).%
\footnote{The  \lq if \rq  and  the \lq only if \rq  parts of Indication property are denoted as  (P2a) and (P2b), respectively, in \cite{halicka2024unified}.}
\end{pozn}
\begin{pozn}
 A weaker version of \ref{P8} where  $(\x_o,\y_o)\succsim (\x_p,\y_p)$ implies $\theta_o^*\le\theta^*_p$, 
 is known in DEA as (weak) monotonicity. This property is met in GS models under mild assumptions that are satisfied by all common path-based models \citep[Remark~8 and Theorem~9]{halicka2024unified}.
\end{pozn}

Let us briefly revisit the reasons for studying this triplet in more detail. The findings of \cite{halicka2024unified} show that none of the three properties are satisfied in the entire class of GS models. Moreover, there are examples of admissible model parameters, where all three properties are met, but it also appears that in most models the three properties fail simultaneously.  

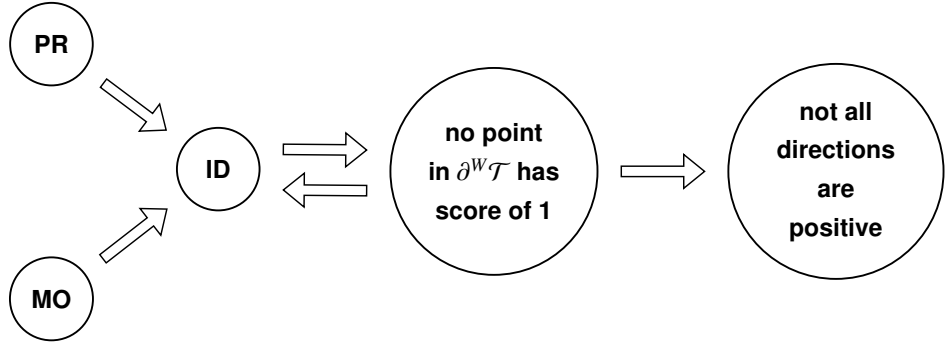
\begin{figure}[ht!]
    \centering
\adjustbox{scale=0.9,center}{\begin{tikzcd}[sep=large,
cells={nodes={circle, draw, thick, inner xsep=6pt}, inner ysep=6pt, outer xsep=10pt},
every arrow/.append style = {-stealth},arrows={tail arrow, tail arrow tip width=1.3em, tail arrow tail width=0.4em, tail arrow tip length=0.8em,tail arrow joint width=0.4em}    ]
\text{\bf \mask{\textsf{PR}}{\textsf{MO}}} \ar[rd]  
&[-1ex]   & [3ex]   &[3ex]   \\[-8ex]
& \text{\bf \mask{\textsf{ID}}{\textsf{MO}}} 
 \arrow[r, yshift=0.3cm] 
&|[circle,inner xsep=2pt, inner ysep = 2pt]| 
\parbox{2.5cm}{\centering \bf\textsf{no point}\\ \textsf{in} $\partial^W\T$ \textsf{has}\\\textsf{score of 1}} \arrow[l, yshift=-0.3cm]\ar[r] 
&|[circle,inner xsep=0pt]| \parbox{2cm}{\centering \bf \textsf{not all}\\\textsf{directions}\\\textsf{are positive}}\\[-8ex]
\text{\bf \textsf{MO}} \ar[ru]  
&   &    &         
\end{tikzcd}
}
    \caption{Hierarchy of properties for 
    the GS model corresponding to a fixed choice of admissible model parameters.}
    \label{Fig:1}
\end{figure}

The relationships among the three properties for a particular GS model with a fixed set of admissible model parameters have been partially explored in \cite{halicka2024unified}. 
Since \ref{P2a} is a restriction of \ref{P3} to the weakly efficient frontier, it is easy to see
that  \ref{P3} $\Rightarrow$ \ref{P2a}, or equivalently, 
 \[\text{NOT \ref{P2a} $\Rightarrow$ NOT \ref{P3}.\footnotemark} \]%
 \footnotetext{By NOT (XY) we understand that the model corresponding to the given  configuration of admissible parameters violates the property (XY) at least at one $(\x_o,\y_o)\in\T$.}%
 NOT \ref{P2a} occurs if and only if there is a unit on the weakly efficient frontier whose score is 1. Hence  \citet[Theorem~4.20]{halicka2024unified} further yields 
 \[\text{NOT \ref{P2a} $\Rightarrow$ NOT \ref{P8}.}\] 
 Moreover, NOT \ref{P2a} arises whenever a unit on the weakly efficient frontier is associated with a positive direction  \citep[Theorem 3]{halicka2024unified}. 
The hierarchy of these properties, formulated by contra-position, is summarised in Figure~\ref{Fig:1}.

It remains to establish the link between (PR) and (MO) and find examples, if any, where (ID) holds but (PR) and/or (MO) do not. This will be the subject of Section~\ref{S3} with the main findings summarised in Figure~\ref{Fig:3}. 

A natural next question is what selection of admissible parameters leads to the success or failure of these characteristics. 
To this end, illustrative examples in  \cite{halicka2024unified} indicate that the success could be connected with a special configuration of the data $\X,\Y$ defining $\T$ together with a special choice of unit-dependent directions based on \citeauthor{portela.al.04}'s \emph{range directions}  
\begin{equation}\label{range directions}
 g^x_o=\x_o- \x^{\m},   \qquad g^y_o=\y^{\M}-\y_o,   
\end{equation} 
suitably modified to fit the GS scheme. This topic is addressed in Section~\ref{S4}.

\section{Connections among \ref{P2a}, \ref{P3}, and \ref{P8}}\label{S3}

Unless the model parameters are explicitly specified, the results of this section apply to a particular GS model with a fixed but arbitrarily chosen set of admissible model parameters. 
We are looking to investigate the consequences of \ref{P2a} on \ref{P3} and \ref{P8}, as well as the connection between \ref{P3} and \ref{P8}. In preparation for these tasks, we recall the necessary and sufficient conditions for \ref{P3}.

   \begin{theorem}[\citealp{halicka2024unified}, Theorem~4]\label{indication}
   The projection of $(\x_o,\y_o)$ is strongly efficient if and only if for each optimal solution $(\la^*_o,\theta^*_o)$ of \eqref{general}, the inequality constraints \eqref{general2} and \eqref{general3} are satisfied with equality.
 \end{theorem}
 
\subsection {\ref{P2a} does not imply \ref{P3}}\label{SS3.1}

The next example shows that \ref{P2a} may hold while \ref{P3} fails for a certain selection
of admissible model parameters.

\begin{figure}[ht!]
    \centering
        \includegraphics[width=10cm]{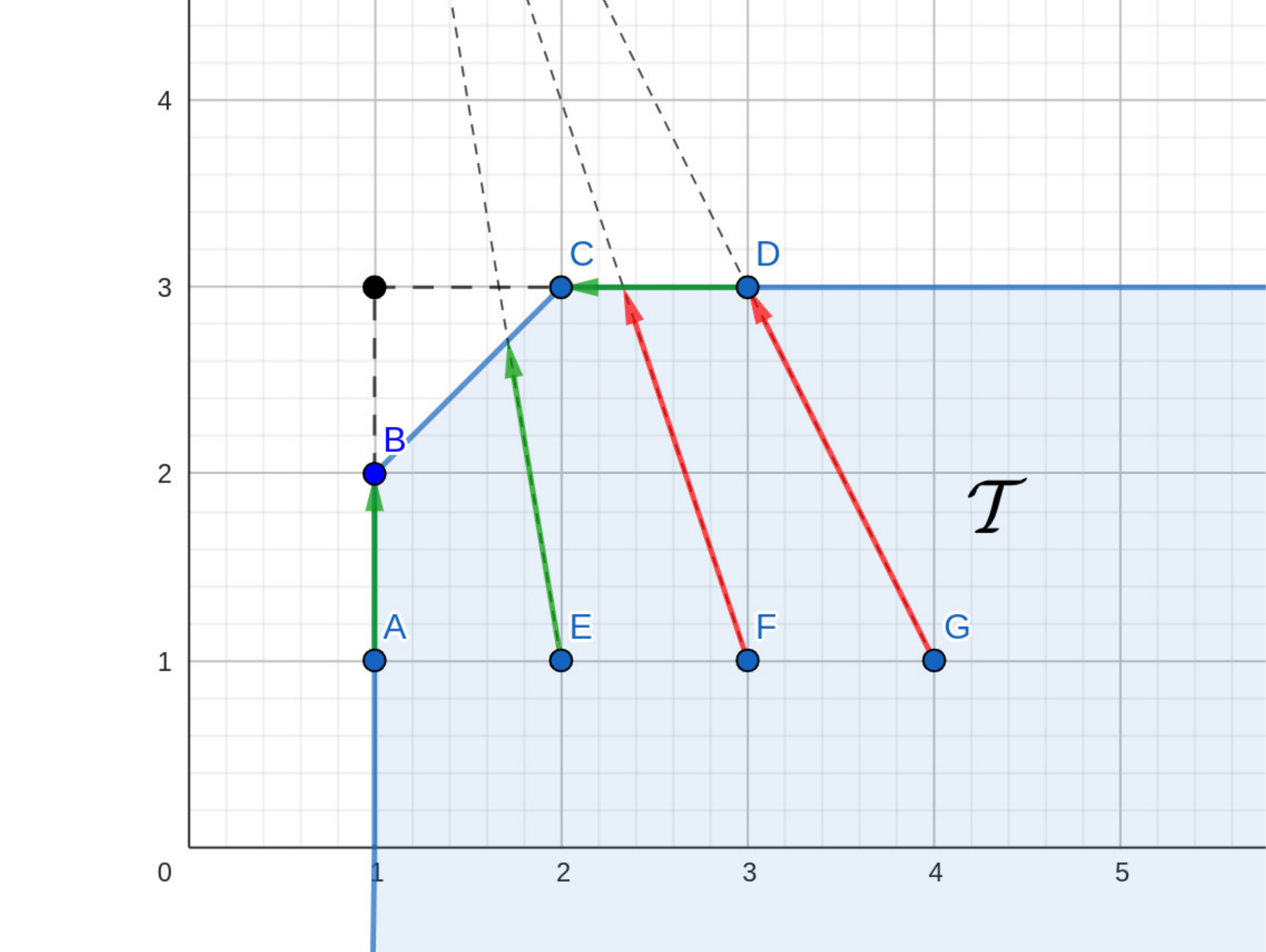}
    \caption{DDF-g model with specific directions over one input and one output technology set. All units belonging to $\partial^W\T$ are projected on $\partial^S\T$, but some units belonging to int$\T$ are projected on $\partial^W\T$. Therefore, \ref{P2a} is satisfied, but \ref{P3} is not.  
    }
    \label{Fig:2}
\end{figure}

\begin{example} \label{ex3.2} Consider a one-input and one-output example with 7 DMUs: units $A=(1,1)$, $B=(1,2)$, $C=(2,3)$, and $D=(3,3)$ on the boundary of $\mathcal{T}$ and units $E=(2,1)$, $F=(3,1)$, and $G=(4,1)$ in its interior. Now, let us apply the DDF-g model with directions  $g^x_o=\x_o- \x^{\m}$,  
$g^y_o=3(\y^{\M}-\y_o)$. Figure~\ref{Fig:2} shows that all units belonging to $\partial^W\T$ (including generic $A,D$) are projected onto $\partial^S\T$,  and the only units that are projected onto themselves, and therefore have a score of 1, are the strongly efficient units from the line segment BC. Therefore, the property \ref{P2a} is satisfied. On the other hand,  some units belonging to $\int\T$ (e.g., $F,G$) are projected onto $\partial^W\T$ and therefore the property \ref{P3} fails.%
\footnote{This example can be adapted to a more general case of $m$ inputs and one output, or one input and $s$ outputs over the so-called `ideal technology sets' introduced later in Section~\ref{S4}.}
Note that a similar situation occurs for the HDF-g model with the \citeauthor{portela.al.04}'s range directions \eqref{range directions}. 
\end{example}

\subsection{\ref{P8} implies \ref{P3}}\label{SS3.2}

Until now, we have considered the right-hand side of
\eqref{general2} and \eqref{general3} as a (vector) function of $\theta$ for a fixed $(\x_o,\y_o)$, which yields a path passing through the assessed unit $(\x_o,\y_o)\in\T$; see \eqref{pathdef}. To obtain a link between \ref{P8} and \ref{P3}, it is helpful to consider the right-hand sides of
\eqref{general2} and \eqref{general3} as functions of the assessed unit $(\x_o,\y_o)$ for a fixed $\theta$. Observe that the right-hand sides in \eqref{pathdef} depend explicitly on $(\x_o,\y_o)$ but the dependence may also be implicit via the directions $\g_o=(\g_o^x,\g_o^y)$.
\begin{definition}\label{D:flow}
 For each fixed $\bar{\theta}\in\cal D$ we define the path-flow mapping $\p( \x_o,\y_o;\bar{\theta})$ as the function that maps each unit $(\x_o,\y_o)\in \T$ to the point $(\p_o^x(\bar{\theta}), (\p_o^y(\bar{\theta}))\in\R^{m+s}$. 
\end{definition}
\begin{theorem}\label{challenge}
    If there exists $(\x_o,\y_o)\in\T$ such that $(\p_o^x(\theta_o^*),\p_o^y(\theta_o^*))\in \partial^W\T$ and the path-flow mapping $\p(\x_o,\y_o;\theta_o^*)$ is continuous at $(\x_o,\y_o)$, then the GS model does not meet the property of strict monotonicity \ref{P8}.
\end{theorem}
\begin{proof}  Since $(\p_o^x(\theta_o^*),\p_o^y(\theta_o^*))\in\partial^W\T$, by Theorem \ref{indication} there exists an optimal solution $(\la^*_o,\theta^*_o)$ of \eqref{general} and $(\bm s^{x*},\bm s^{y*})\ge \0$ such that $(\bm s^{x*},\bm s^{y*})\neq \0$ and
\begin{equation}\label{phase1}
\X\la^* +\s^{x*} =\p_o^x(\theta^*_o), \quad \Y\la^* -\s^{y*} =\p_o^y(\theta^*_o).\\
\end{equation}
Define 
$${\cal S}=\{(\x_o-\s^x, \y_o+\s^y) \ \vert \  0\lneqq (\s^x,\s^y)\le (\s^{x*},\s^{y*}) \}.$$
Each unit in $\cal S$ is distinct from and dominates $(\x_o,\y_o)$. Since $(\x_o,\y_o)=\p_o(1)$ and $\theta^*_o\le 1$, the monotonicity of the path in $\theta$ implies that $\p_o^x(\theta^*_o)\le \x_o$ and $\p_o^y(\theta^*_o)\ge \y_o$. Therefore, the equations in \eqref{phase1} imply that
$(\x_o-s^{x*}, \y_o+s^{y*})\in \T $ and therefore also ${\cal S}\subset \T$. Note that $\cal S$ contains points that are arbitrarily close to $(\x_o,\y_o)$. Now, by the continuity assumption, for $\epsilon:=\min_{i,r} \{s^{x*}_i>0, s^{y*}_r>0\}$ there exists a unit $(\x_p,\y_p)\in\cal S$ such that $\vert\vert\p_o(\theta^*_o)-\p_q(\theta^*_o)\vert\vert<\epsilon$. This implies that the vector inequalities $\vert\p^x_o(\theta^*_o)-\p^x_p(\theta^*_o)\vert\le\s^{x*}$ and $\vert\p^y_o(\theta^*_o)-\p^y_p(\theta^*_o)\vert\le\s^{y*}$ hold.  Therefore, one also has
\begin{equation}\label{phase3}
\p_o^x(\theta^*_o)-\s^{x*}\le \p_p^x(\theta^*_o),\quad \p_o(\theta^*_o)+\s^{y*}\ge \p_p^y(\theta^*_o).
\end{equation}
By substituting \eqref{phase3} into \eqref{phase1} we get
\begin{equation}\label{phase2}
      \X\la^* \le\p_p^x(\theta^*_o), \quad
      \Y\la^*  \ge\p^y_p(\theta^*_o).
\end{equation}
This yields that $(\la^*,\theta_o^*) $ is a feasible solution for (GS)$_p$, and hence one has $\theta_p^*\le\theta_o^*$.
The strict monotonicity is violated because $(\x_p,\y_p)$ dominates and is different from $(\x_o, \y_o)$ but $\theta^*_p\le\theta_o^*$.
\end{proof}

Observe that Theorem~\ref{challenge} formulates a  local property: if a point in $\T$ is \emph{not} projected onto $\partial^S\T$, then one can find another point such that the pair fails to maintain strict monotonicity. Thus, the ``local'' formulation offers more flexibility than the following ``global'' corollary.

\begin{corollary}\label{C:5.3} Assume that the path-flow mapping $\p( \x_o, \y_o; \theta)$ is continuous  for each fixed $\theta\in\D$ on $\T$. Then \ref{P8} $\Rightarrow$ \ref{P3} holds.
\end{corollary}

\subsection{\ref{P3} implies \ref{P8}}\label{SS3.3}

The monotonicity of GS models is related to the monotonicity of the path-flow mapping $\p(\x_o,\y_o; \theta)$, which we formalise next.
\begin{definition}\label{defmonoton}
We say that the path-flow mapping $\p( \x_o,\y_o;\bar{\theta})$: $(\x_o,\y_o)\to \p_o(\bar{\theta})$ is \emph{monotone} on $\T$ at $\bar{\theta}\in \mathcal{D}$ if for any two units $(\x_o,\y_o)$, $(\x_p,\y_p)$ in $\T$, one has
\begin{equation}\label{monoton}
    (\x_o,\y_o)\succsim (\x_p,\y_p) \ \Rightarrow
    \ \p_o(\bar{\theta})\succsim \p_p(\bar{\theta}).
\end{equation}
If, in addition,  
 \begin{equation}\label{smonoton}
    (\x_o,\y_o)\succnsim (\x_p,\y_p) \ \Rightarrow
   \ \p_o(\bar{\theta})\succnsim \p_p(\bar{\theta}),
\end{equation}
we  say that the path-flow mapping is \emph{strictly monotone} on $\T$ at $\bar{\theta}\in \mathcal{D}$.
\end{definition}
\begin{pozn}\label{pozlozkach} If for each $i\in\{1,\dots,m\}$ the 
$i$-th component $[\p^x( \x_o,\y_o;\bar{\theta})]_i$ of $\p^x( \x_o,\y_o;\bar{\theta})$
depends only on $x_{io}$, and for each $r\in\{1,\dots,s\}$ component $[\p^y( \x_o,\y_o;\bar{\theta})]_r$ depends only on
$y_{ro}$, then the
monotonicity property in Definition~\ref{defmonoton}
simply means that $[\p^x(\bar{ \x_o,\y_o;\theta})]_i$ and  $[\p^y( \x_o,\y_o;\bar{\theta})]_r$ are nondecreasing
in $x_{io}$ and $y_{ro}$, respectively. On the other hand, the strict monotonicity property in Definition~\ref{defmonoton}
 means that $[\p^x( \x_o,\y_o;\bar{\theta})]_i$ and $[\p^y( \x_o,\y_o;\bar{\theta})]_r$ are increasing
in $\x_{io}$ and $\y_{ro}$, respectively.
\end{pozn}

\begin{pozn}\label{poznP7} Note that if  $\g_o$
does not depend on  $(\x_o,\y_o)$, then the path-flow mapping $\p( \x_o,\y_o;\theta)$ is monotone on $\T$ for any choice of $\psi^x$ and $\psi^y$ at any $\theta\in \mathcal{D}$.
\end{pozn}
The next lemma will be useful for further analysis of strict monotonicity.

\begin{lemma} [\citealp{halicka2024unified}, Lemma~3]\label{monotonlema}
Let $(\x_o,\y_o)$ and $(\x_p,\y_p)$ be two units in $\T$
with the corresponding optimal values $\theta^*_o$ and $\theta^*_p$. If $\p_o(\theta^*_o)\succsim \p_p(\theta^*_o)$, then each optimal solution
$(\theta_o^*,\la_o^*)$  of (GS)$_o$ is a feasible solution of (GS)$_p$ and hence $\theta^*_p\le \theta^*_o$.
\end{lemma}
Let us note that the assumption of monotonicity of the path-flow mapping ensures the property of the so-called weak monotonicity of each model in the GS scheme as indicated by Lemma \ref{monotonlema}. To establish the strict monotonicity \ref{P8} we will need the strict monotonicity of the path-flow mapping and property \ref{P3}. First, we formulate a local version of the assertion.
\begin{theorem}\label{strongmonotonicity2}
Suppose that units $(\x_o,\y_o)$  and $(\x_p,\y_p)$ in $\T$ with the efficiency scores $\theta^*_o$ and $\theta^*_p$, respectively, are projected onto the strongly efficient frontier. 
Suppose also that $(\x_o,\y_o)\succnsim (\x_p,\y_p)$ and that the path-flow mapping is strictly monotone at $\theta_o^*$. Then $\theta_o^*>\theta^*_p$. 
\end{theorem}
\begin{proof}
 From the assumptions of the theorem it follows that $\p_o(\theta_o^*)\succnsim \p_p(\theta_o^*)$. This, by Lemma \ref{monotonlema} implies that $\theta^*_o\ge \theta^*_p$. Assume, by contradiction, that $\theta^*_p= \theta^*_o$. Denote by $(\la_o^*,\theta_o^*)$  an optimal solution for (GS)$_o$. Obviously, the same pair also represents an optimal solution for (GS)$_p$. Since the model projects $(\x_o,\y_o)$ and $(\x_p,\y_p)$ onto the strongly efficient frontier, by Theorem \ref{indication} the inequalities \eqref{general2} and \eqref{general3} in (GS)$_o$ and (GS)$_p$ are binding, that is, 
  \begin{alignat}{3}
  &\X\la_o^*&{}=\p_o^x(\theta_o^*) &\text{\quad and\quad }&  \Y\la_o^*&{}=\p^y_o(\theta^*_o).\label{1}
  \end{alignat}
    \begin{equation}\label{3}
   \X\la_o^*=\p_p^x(\theta_o^*)\text{\quad and \quad}  \Y\la_o^*=\p_p^y(\theta^*_o).
\  \end{equation}
By comparing the right-hand sides of \eqref{1} and \eqref{3},
we get
   \begin{equation}\label{4}  \p_o^x(\theta_o^*)=\p^x_p
  (\theta_o^*), \quad
  \p^y_o(\theta^*_o)=\p_p^y
   (\theta^*_o).
  \end{equation}
  This is in contradiction with the assumption of strict monotonicity of the path-flow mapping at $\theta_o^*$, according to which $\p_o(\theta_o^*)\succnsim \p_p(\theta_o^*)$ holds.
\end{proof}
Theorem~\ref{strongmonotonicity2}, too, is formulated locally, i.e., if two  ordered units are both projected onto $\partial^S\T$, then their efficiency scores are strictly ordered. The local formulation once again offers more flexibility than the following global corollary. 

\begin{figure}[ht!]
 \centering
    \adjustbox{scale=0.9,center}{\begin{tikzcd}[cells={nodes={circle, draw, thick, inner xsep=6pt}, inner ysep=6pt,outer xsep=14pt}, row sep={0cm,between origins},column sep={0cm,between origins},arrows={tail arrow, tail arrow tip width=1.3em, tail arrow tail width=0.4em, tail arrow tip length=0.8em,tail arrow joint width=0.4em}]
\text{\bf \textsf{MO}} \ar[r, shift left=0.3cm, "{\raisebox{1.1em}{\text{\bf \textsf{path-flow continuity}}}}"] 
&[6cm] \text{\bf \mask{\textsf{PR}}{\textsf{MO}}} \ar[l, shift left=0.3cm, "{\raisebox{-1.8em}{\text{\bf \textsf{path-flow monotonicity}}}}"] \ar[r] &[5cm]  \text{\bf \mask{\textsf{ID}}{\textsf{MO}}}
\end{tikzcd}
}
\caption{Hierarchy of three properties in the models of the GS scheme.}
\label{Fig:3}
\end{figure}
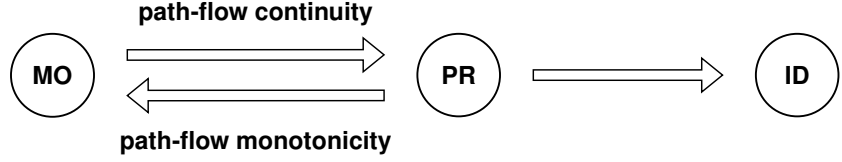

 \begin{corollary} \label{P3P8} Denote by ${\cal D}^*$ the set of efficiency scores
 $\theta_o^*$ achievable by units 
$(\x_o,\y_o)\in\T$. Assume that the path-flow mapping $\p( \x_o,\y_o;\theta^*)$ is strictly monotone on $\T$ for all $\theta^*\in{\cal D}^*$.  Then \ref{P3} $\Rightarrow$ \ref{P8} holds.
\end{corollary}
The results of this section are summarised in Figure~\ref{Fig:3}.

\section{Specific technology sets and direction vectors ensuring  \ref{P3}}\label{S4}

We recall that the GS models fail all three properties \ref{P2a}, \ref{P3}, and \ref{P8} when the directions are positive. Nonetheless, there exist simple examples, where 
the GS models with specific directions over specific technology sets project every unit from $\T$ onto the strongly efficient frontier.
\begin{figure}[htbp] \label{dvojobrazok}
    \centering
        \includegraphics[width=7.5cm]{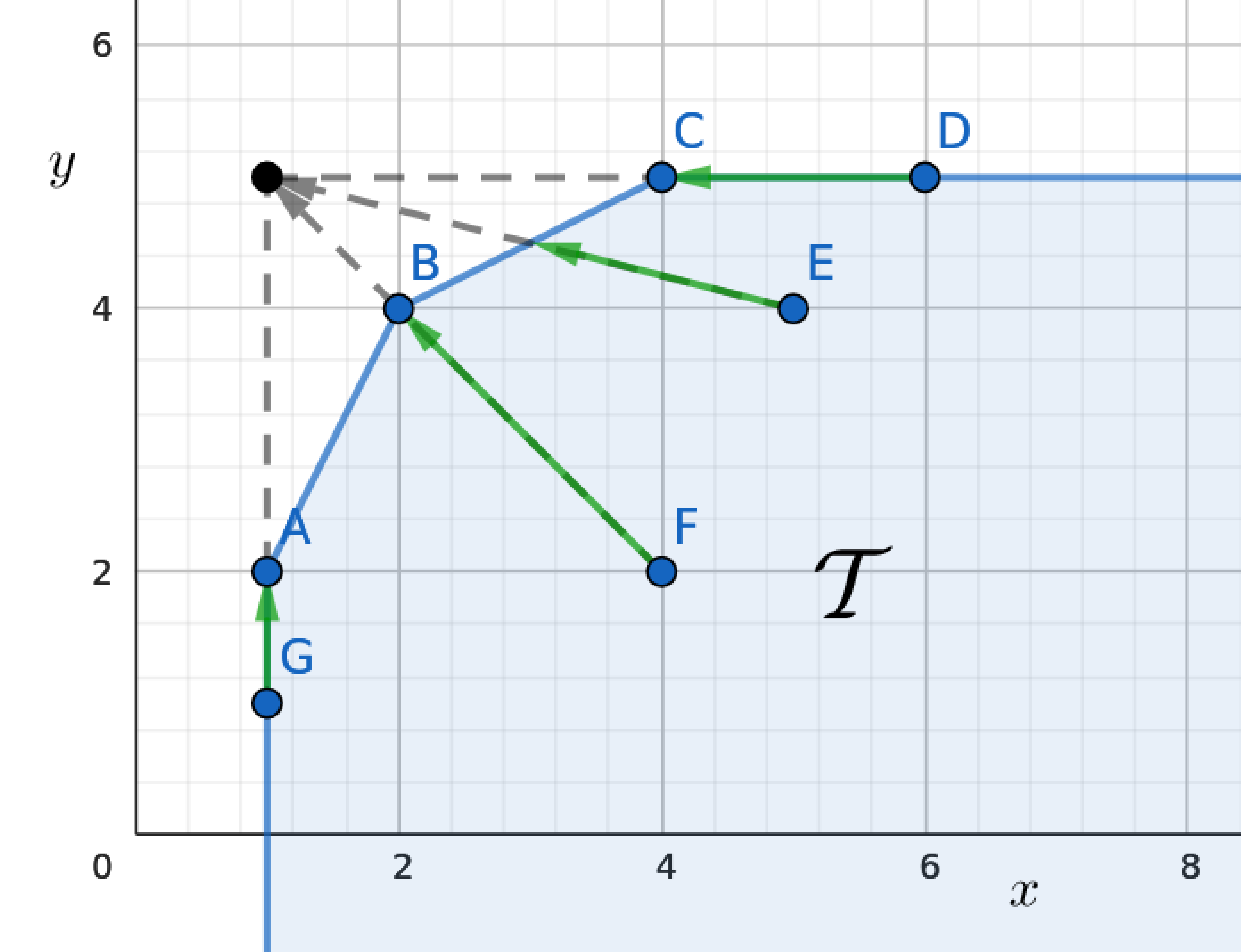}
         \includegraphics[width=7.5cm]{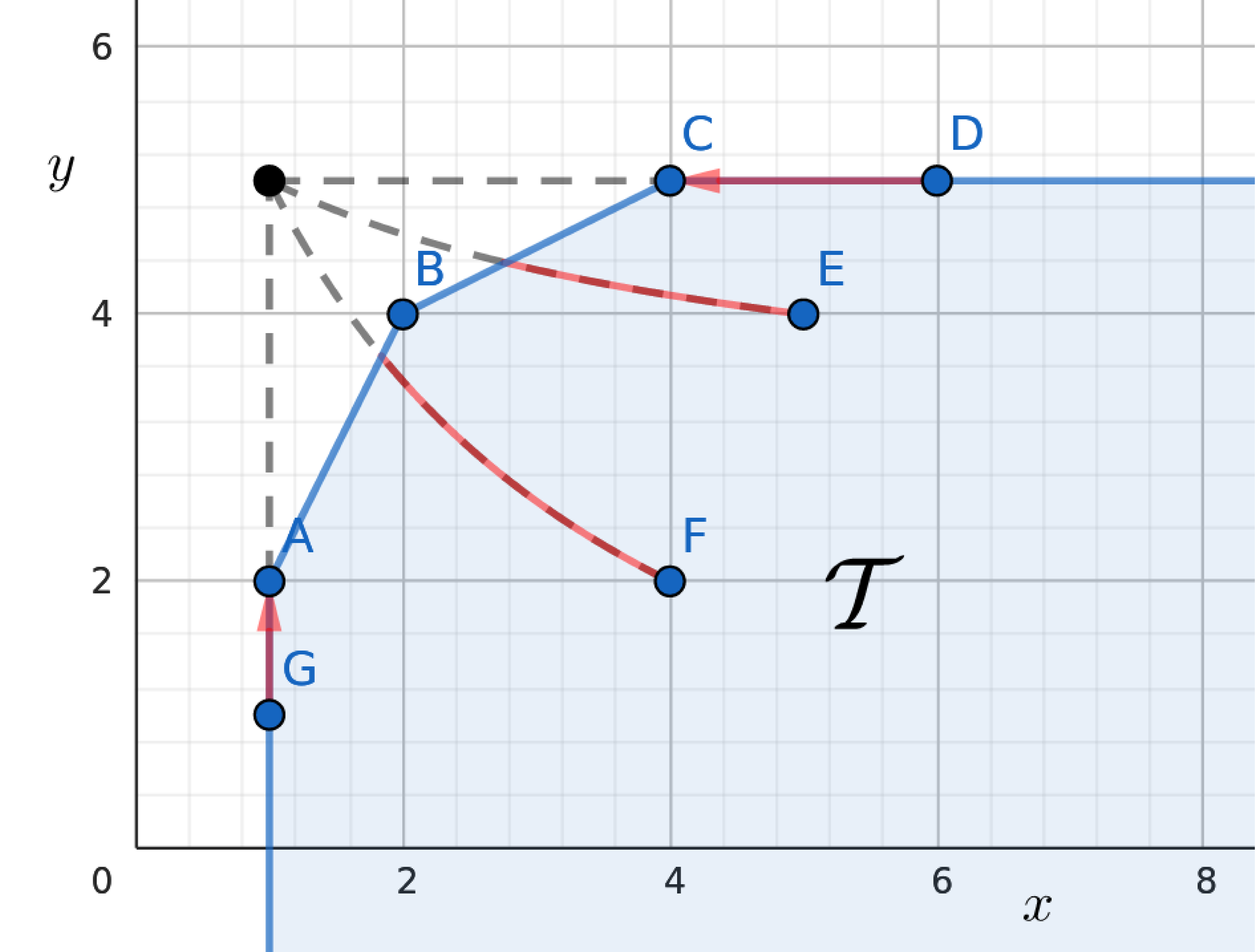}
    \caption{Illustration of Example \ref {E4.1}. GS models in the case of a technology set $\T$ with one input and one output that project all units on $\partial^S\T$. Left: DDF-g paths with the range directions \eqref{range directions}. Right: HDF-g paths with the modified range directions. In both cases, the paths pass through the ideal point $(\x^\m,\y^\M)$.}
    \label{Fig:4}
\end{figure}
\begin{example}\label{E4.1} Consider a one-input, one-output technology set generated by three units $A,B,C$ as depicted in Figure~\ref{dvojobrazok}.  DDF-g paths (where $\psi^x(\theta)=\theta$ and $\psi^y(\theta)=2-\theta$) with the range directions $\g^x_o=\x_o-\x^{\m}$ and $\g^y_o=\y^{\M}-y_o$ connect the units from $\T$ with the ideal point $(1,5)$ for this technology as seen in the left diagram of Figure \ref{dvojobrazok}. In the case of  HDF-g paths (where $\psi^x(\theta)=\theta$ and $\psi^y(\theta)=\frac{1}{\theta}$) the range directions were modified to $\g^x_o=2(\x_o-\x^{\m})$ and $\g^y_o=\y^{\M}-y_o$ to pass through the ideal point $(1,5)$ as shown in the right panel of  Figure~\ref{dvojobrazok}. In both cases, the models project units onto the strongly efficient boundary consisting of the line segments $AB$ and $BC$.
\end{example}

 In this section, we will generalise the two GS models presented in Example \ref{E4.1} with the aim of fully characterising those technology sets and directions under which the GS model yields the strong efficiency of projections. In doing so, the ideal point $(\x^\m,\y^\M)$ of the technology set  $\T$  will play an important role. 

\subsection{GS range directions}\label{SS4.2}
 
As shown in Figure~\ref{Fig:4}, the paths in Example~\ref{E4.1} connect units from $\T$ to the ideal point $(\x^\m,\y^\M)$, which ensures that the units are projected onto the strongly efficient boundary. Inspired by this example,
we now characterise those directions $g_o$ for which the corresponding path $\p_o$ for  $(\x_o,\y_o)$   passes through the ideal point $(\x^\m,\y^\M)$ at an arbitrarily chosen $\theta_{\m}\in [0,1)$. The choice of common $\theta_{\m}$ for all units in $\T$ then allows us to obtain comparable scores whose common lower bound is $\theta_\m$.  The proof follows by a simple calculation and is therefore omitted. 

Recall that $(\x^\m,\y^\M)$ denotes the ideal point of $\T$ and that for fixed $\psi^x$ and  $\psi^y$,  the path $\phi_o$ is determined by the evaluated unit $(x_o,y_o)$ and the direction $\g_o=(\g_o^x, \g_o^y)$ as detailed in Subsection~\ref{SS2.2}.

\begin{lemma}\label{dir} Fix $\psi^x$, $\psi^y$ satisfying assumptions \ref{A1}--\ref{A3} 
and $\T$ of the form \eqref{T}. For $\theta_{\m}\in [0,1)\cap\D $ and $(\x_o, \y_o)\in\T\setminus\{ (\x^\m, \y^\M)\}$, the following are equivalent.
\begin{enumerate} [(i)]
\item The path $\p_o$ runs through the ideal point $(\x^\m, \y^\M)$ at $\theta=\theta_{\m}$, i.e., $\p_o({\theta}_{\m})=(\x^\m,\y^\M)$.
\item The direction vector $\g_o=(\g_o^x, \g_o^y)$ of $\p_o$ satisfies \begin{equation}\label{directions}
 \g_o^x=\frac{\x_o-\x^\m}{1-\psi^x({\theta}_{\m})}, \quad
 \g_o^y=\frac{\y^\M-\y_o}{\psi^y({\theta}_{\m})-1}.
 \end{equation}
\end{enumerate}

 Moreover, if either of these conditions holds, then $\theta_{\m}\le\theta^*_o$.
\end{lemma}
We will call the direction vector defined by \eqref{directions} the \emph{GS range direction}.

In the case of DDF-g models, we have $\psi^x(\theta)=\theta$, $\psi^y(\theta)=2-\theta$, and $\mathcal D=\R$. On selecting $\theta_{\m}=0$, the GS range directions coincide with the range directions \eqref{range directions} of \cite{portela.al.04}. In the case of HDF-g models, we have $\psi^x(\theta)=\theta$, $\psi^y(\theta)=\frac{1}{\theta}$, and $\mathcal D=(0,\infty)$, hence $\theta_{\m}$ must be chosen positive. The value $\theta_{\m}=1/2$ yields $\g^x_o=2(\x_o-\x_{\m})$ and $\g^y_o=\y^{\M}-\y_o$.

\begin{pozn}
For $(\x_o, \y_o)\in\T\setminus\{ (\x^\m, \y^\M)\}$, the GS range  direction for $(\x_o, \y_o)$ satisfies $\g_o=(\g_o^x,\g_o^y)\gneqq 0$ and therefore the programme $(GS)_o$ is well defined.  On the other hand, if $(\x^\m, \y^\M)\in\T$, then the GS range direction for $ (\x^\m, \y^\M)=(\x_o,\y_o)$ vanishes and the corresponding $(GS)_o$ programme is not well defined. However, since $ (\x^\m, \y^\M)=(\x_o,\y_o)$ is the (only) strongly efficient unit in $\T$, we can set $\theta_o^*=1$ by definition for this point. 
\end{pozn}
The next lemma shows that the choice of the GS range direction ensures  that both $\p_o^x$ and $\p_o^y$, the $x$ and $y$ components of the path  $\p_o$, passes through the relative interior of the  line segments connecting  $\x_o$ with $\x^{\m}$ and $\y_o$ with $\y^{\M} $, respectively.

\begin{lemma}\label{dirmod} Let $\theta_{\m}\in [0,1)\cap\D $, $(\x^\m, \y^\M)\notin \T$ and  $\theta_o^*< 1$ be the optimal value of $(\x_o, \y_o)\in\T$ in the $(GS)_o$ model 
 with the GS range direction. Then, for each  $\hat\theta\in[\theta_o^*,1)$ 
  there exist $\alpha^x(\hat\theta)\in(0,1)$ and $\alpha^y(\hat\theta)\in(0,1)$ such that
 \begin{equation}\label{pathmod}
 \p_o^x(\hat\theta)=(1-\alpha^x(\hat\theta))\x_o+\alpha^x(\hat\theta)\x^m, \quad
 \p_o^y(\hat\theta)=(1-\alpha^y(\hat\theta))\y_o+\alpha^y(\hat\theta)\y^M.
 \end{equation} 
\end{lemma}
\begin{proof}
For each $\hat\theta\in[\theta_o^*,1)$ we have $(\x^{\m},\y^{\M})=\p_o(\theta_{\m})\succnsim \p_o(\theta^*)\succeq\p_o(\hat\theta)\succnsim \p_o(1)=(\x_o,\y_o)$ and hence $\psi^x(\theta_{\m})<\psi^x(\theta^*_o)\le \psi^x(\hat\theta)< 1$ and $\psi^y(\theta_{\m})>\psi^y(\theta^*_o)\ge\psi^y(\hat\theta)> 1$. 
 Therefore,
 \begin{equation}\label{alfa}
\alpha^x(\hat\theta):=\frac{1-\psi^x(\hat\theta)}{1-\psi^x({\theta}_{\m})}\in(0,1), \quad
\alpha^y(\hat\theta):=\frac{\psi^y(\hat\theta)-1}{\psi^y(\theta_{\m})-1}\in(0,1).
\end{equation} 
 Finally, a simple computation yields 
\begin{equation*}
\begin{aligned}
 \p_o^x(\theta^*_o)=\x_o-\alpha^x(\hat\theta)\g_o^x=(1-\alpha^x(\hat\theta))\x_o+\alpha^x(\hat\theta)\x^{\m}, \\
 \p_o^y(\theta^*_o)=\y_o+\alpha^y(\hat\theta)\g_o^y=(1-\alpha^y(\hat\theta))\y_o+\alpha^y(\hat\theta)\y^{\M},
  \end{aligned}
 \end{equation*}
 from which the claim of the lemma follows.
\end{proof}

Figure \ref{dvojobrazok} illustrates that the paths containing the ideal point of a given technology set project all units onto the strongly efficient frontier. This property, whilst true in  two-dimensional technology sets, no longer holds for all higher-dimensional technologies. Indeed, the two-input, one-output example in \citet[Example~2 and Figure 3]{halicka2024unified} shows that the DDF-g model with the GS range directions projects a unit from $\partial^W\T$ onto itself. 
\cite{halicka2024unified} conducted further numerical experiments on real data consisting of 30 units, 2 inputs, and 2 outputs, with some negative output values. Their analysis contrasts the so-called proportional directions $\g_o^x=\vert\x_o\vert$, $\g_o^y=\vert\y_o\vert$ of \cite{kerstens2011negative} with the GS range directions for four different choices of the path function $\psi$.  In their analysis, the GS range directions yield a significantly smaller proportion of units \emph{not} projected onto the strongly efficient frontier but do not eliminate such units completely. 
Therefore, even though the GS range directions ensure the passage of all paths through the ideal point, this alone is not enough to guarantee that units are projected onto the strongly efficient frontier in the case of higher-dimensional technology sets. 
\subsection{Ideal technology sets}\label{ss42}

It is obvious that to ensure the fulfilment of (PR) using GS range directions, it will be necessary to limit ourselves to a certain subclass of technology sets. To specify the properties of the subclass, it turned out that the facial structure of the technology as a polyhedral set will play an important role. Some characteristics of the facial structure of \eqref{1}, as well as all necessary terms, are included in Appendix A. The results of the current section will finally show that the GS range directions will guarantee projecting onto strongly efficient frontier if and only if we restrict ourselves to technology sets specified in the following definition.  
\begin{definition}\label{D:idealT}
A technology set $\T$ of the form \eqref{T} is called an \emph{ideal technology} if each unit $(\x,\y)$ from the weakly efficient boundary of $\T$    
has at least one component in common with the ideal point $(\x^\m,\y^\M)$ of  $\T$, i.e.,
\begin{equation}\label{PW}
\forall\ (\x,\y)\in \partial^W \T: \ (\exists \ i\in \{1,...,m\}: \ x_i=x_i^\m) \ \hbox{or} \
(\exists\ r\in \{1,...,s\}: \ y_r=y_r^\M). 
\end{equation}
\end{definition}

It is easy to see that the technology sets \eqref{T} with only one input and one output, as in Example~\ref{E4.1}, are ideal. A trivial example of an ideal technology in any dimension is the trivial technology, i.e.,  the technology containing its ideal point $(\x^\m,\y^\M)$. The trivial technology is a polyhedral cone with the vertex at the ideal point. Its facets --- faces of full dimension --- are parallel to the orthant hyperplanes. The vertex is the only strongly efficient unit of the trivial technology set. 

\begin{figure}[htb]
    \centering
       {\includegraphics[width=7cm]{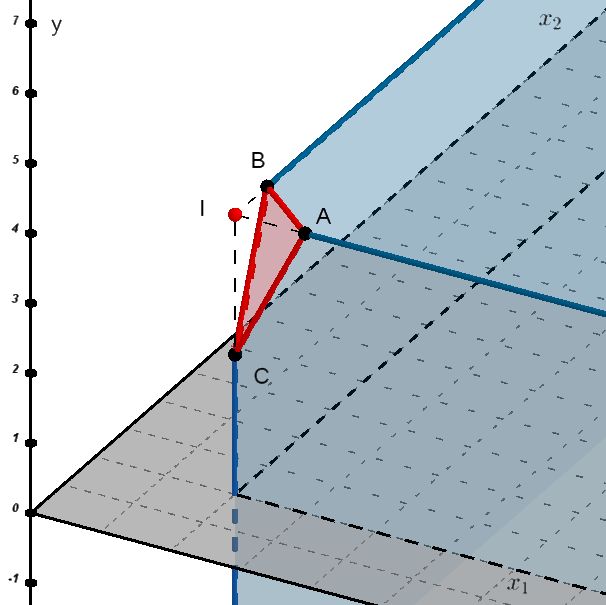}}
    \caption{The ideal technology set from Example \ref{tech1}. It contains just $m+s=3$  unbounded (blue) edges. }
    \label{Fig:5}
\end{figure}

 We now provide examples of ideal and non-ideal technology sets.
\begin{example} \label{tech1} Consider two inputs and one output example with three DMU`s: $A=(3,2,4)$, $B=(2,3,4)$, $C=(2,2,2)$, see Figure~\ref{Fig:5}. 
All three units $A, B$ and $C$ are strongly efficient and the triangle
$ABC$ forms a strongly efficient frontier for the technology generated by these three DMU's. Apparently
$I=(\x^\m, \y^\M)=(2,2,4)$. It is seen that the closure of $\partial^W\T$ consists of three facets, each determined by one of the hyperplanes $x_1=2$, $x_2=2$ or $y=4$, and thus this technology set satisfies the property \eqref{PW}.
\end{example}

\begin{figure}[htb]
    \centering
        {\includegraphics[width=8cm]{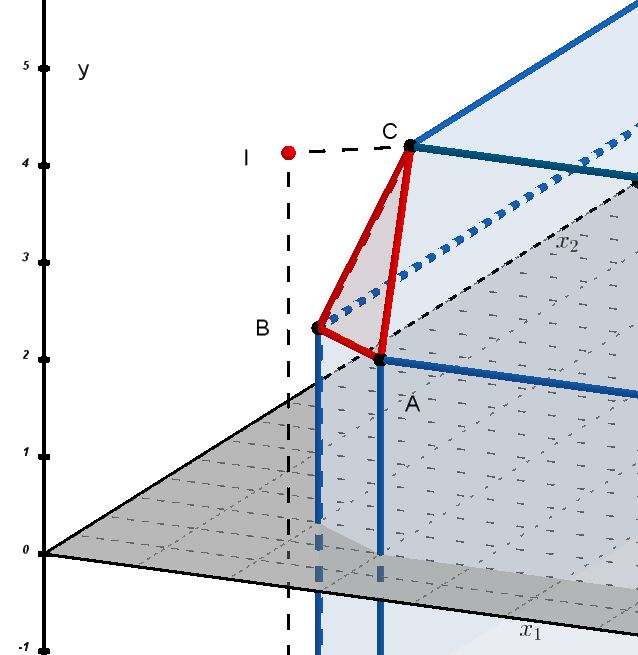}}
    \caption{The non ideal technology set from Example~\ref{tech2}. It contains six unbounded (blue) edges.}
    \label{Fig:6}
\end{figure}

\begin{example}\label{tech2}
Consider again two inputs and one output example with another three DMU`s: $A=(3,2,2)$, $B=(2,3,2)$, $C=(3,3,4)$, see Figure~\ref{Fig:6}. 
Also in this example all three points are strongly efficient and triangle $ABC$ forms a strongly efficient frontier for the technology generated these three DMU's. Apparently $I=(\x^\m, \y^\M)=(2,2,4)$. In this case, the closure of $\partial^W\T$ consists of six facets, where only three of them correspond to hyperplanes $x_1=2$, $x_2=2$ or $y=4$.   Points $(x_1,x_2,y)$ from (the relative) interiors of the other three facets satisfy $x_1>x_1^{\m}=2$, $x_2>x_2^{\m}=2$ and $y<y^{\M}=4$ and hence property \eqref{PW} is not satisfied.\end{example}

Theorem~\ref{alternativy} provides equivalent characterisations of ideal technologies. 
It shows that every ideal technology is obtained by `tapering off' some part of the trivial technology near the ideal point by means of a finite number of hyperplanes in such a way that all $m+s$  edges of the trivial technology are preserved sufficiently far from its vertex (Figure~\ref{Fig:5}). A technology is not ideal precisely when at least one of the $m+s$ edges of the trivial technology is completely missing. Figure~\ref{Fig:5} gives an example of a non-ideal technology, where all three edges of the induced trivial technology are absent.
Item~\ref{ek.v} of Theorem~\ref{alternativy} also provides a tool for the practical recognition of whether a given set of DMUs generates an ideal technology or not. 

The following lemma is an important instrument for developing the main results in Subsection~\ref{SS4.3}. It characterises the relative boundary of unbounded faces of an ideal technology $\T$. Its proof is found in \ref{C}. The necessary material describing the facial structure of $\T$ as a polyhedral set is collected in \ref{A} and \ref{B}. 

\begin{lemma}\label{AcorPW}
Let $\T$ be an ideal technology set. For $(\hat{\x},\hat{\y})\in \mathcal{T}$  we introduce index sets $I=\{i:\hat x_{i}=x_i^\m\}$ and $R=\{r:\hat y_{r}=y_r^\m\}$ and denote their cardinality by $\vert I\vert$ and $\vert R\vert$, respectively. Then  $(\hat{\x},\hat{\y})$  belongs either to $\partial^S\mathcal{T}$, or to the relative interior of the $(m+s-\vert I\vert-\vert R\vert)$ dimensional face of $\T$ given by 
\begin{equation}\label{PIR}
 \mathcal \F_{IR}:=\{(x,y)\in \T: x_{i}=x_i^\m, i\in I;\ y_{r}=y_r^\m, r\in R\} .  
\end{equation}
\end{lemma}

\subsection{GS models with GS range directions over ideal technology sets} \label{SS4.3}

The next theorem shows that the GS models with GS range directions meet the property \ref{P3} for any feasible $\psi$ if and only if $\T$ is ideal. 

\begin{theorem}\label{T:strongly}
In the GS model with the GS range directions \eqref{directions}, the following are equivalent.
\begin{enumerate}[(i)]
\item\label{T:strongly.i} All units in $\T$ are projected onto the strongly efficient frontier.
\item\label{T:strongly.ii} $\T$ is an ideal technology set.
\end{enumerate}
\end{theorem}

\begin{proof} 
\ref{T:strongly.i} $\Rightarrow$ \ref{T:strongly.ii} Arguing by contradiction, assume that $\T$ is not an ideal technology. This means that $\T$ does not have the property \eqref{PW} and therefore there exists $(\x_o,\y_o)\in\partial^W\T$ such that $\x_o>\x^\m$ and $\y_o<\y^\M$. Then the formulas in \eqref{directions} imply $g_o>0$. Now, by Theorem 4.2 of \cite{halicka2024unified}, the positivity of the direction at  $(\x_o,\y_o)\in\partial^W\T$ implies  $\theta_o^*=1$. Hence $(\x_o,\y_o)\in\partial^W\T$ is identical to its projection and \ref{P3} is not satisfied.  \\[0.5ex] 
\noindent \ref{T:strongly.ii} $\Rightarrow$ \ref{T:strongly.i} If $(\x^\m, \y^\M)\in \T$ (that is, if $\T$ is the trivial technology), then the choice of the GS range direction guarantees that each $(\x_o, \y_o)\in\T\setminus\{ (\x^\m, \y^\M)\}$ is projected on $(\x^\m, \y^\M)\in \T$, which is the only strongly efficient unit in $\T$, and hence the theorem holds. Obviously $\theta^*_o=\theta_{min}$.

  Now consider the case $(\x^\m, \y^\M)\notin \T$.
   Let $I$ and $R$ be the index sets of $({\x_o},{\y_o})\in \mathcal{T}$ defined in Lemma \ref{AcorPW}, which can also be empty. According to Lemma \ref{AcorPW}, $(\x_o,\y_o)$ either belongs to the relative interior of $\mathcal F_{IR}$ or is strongly efficient.  In the latter case, there is nothing to prove.  Therefore, assume that $(\x_o,\y_o)\in \text{rel int}\mathcal F_{IR}$ and $\theta_o^*$ is its score.   From the prescription \eqref{directions} of the GS range direction it follows that $\g_{io}^x=0$ and $\g_{ro}^y=0$  if and only if $i\in I$ and $r\in R$ respectively. This implies that the path $(\p_o^x(\theta), \p_o^y(\theta))$ for decreasing values of $\theta\le 1$ stays in $\mathcal F_{IR}$ until it reaches the relative boundary of $\mathcal F_{IR}$ at some $\hat\theta<1$. This allows us to apply Lemma \ref{dirmod} according to which there exist $\alpha^x(\hat\theta)\in (0,1)$ and $\alpha^y(\hat\theta)\in (0,1)$  such that   $\p_o^x(\hat\theta)$ and  $\p_o^y(\hat\theta)$ satisfy \eqref{pathmod}. From this follows that  $\phi_{io}^x(\hat\theta) =x_i^{\m}$ and  $\phi_{io}^y(\hat\theta) =y_i^{\M}$ if an only if $i\in I$ and $r\in R$, resp.,
  hence $I$ and $R$ is the index set also for $(\p_o^x(\hat\theta),\p_o^y(\hat\theta))$. We apply Lemma \ref{AcorPW} again, this time to a point $(\p_o^x(\hat\theta),\p_o^y(\hat\theta))$ that we know is from the relative boundary of $\mathcal F_{IR}$ and hence $(\p_o^x(\hat\theta),\p_o^y(\hat\theta))$ is strongly efficient and $\theta_o^*=\hat \theta$.
 \end{proof}
\begin{theorem}\label{strictly} Let $\T$ be an ideal technology set with $(\x^\m, \y^\M)\notin \T$. Then the GS model with the GS range direction meets the property of strict monotonicity \ref{P8}.
\end{theorem}
\begin{proof} According to Theorem \ref{T:strongly}, the GS model with the GS range direction satisfies \ref{P3}, and hence by Corollary \ref{P3P8}  it suffices to prove that the path-flow mapping $\p(\theta^*)$ is strictly monotone on $\T$ at any $\theta^*\in{\cal D}^*$. In our case, the vector function $\p({\theta^*})$ has a property that  $[\p^x({\theta^*})]_i$
depends only on $x_{io}$ and $[\p^y({\theta^*})]_r$ depends only on
$y_{ro}$, and hence by Remark \ref{pozlozkach} it suffices to prove that $[\p^x({\theta^*})]_i$ and $[\p^y({\theta^*})]_r$ are increasing
in $x_{io}$ and $y_{ro}$, respectively. Obviously, the path-flow mapping $\p(\theta^*)$ is strictly monotone on $\T$ at  $\theta^*=1$. For the case $\theta^*_o<1$ we use  Lemma \ref{dirmod} according which there exist $\alpha^x\in (0,1)$ and $\alpha^y\in (0,1)$ such that
\begin{equation}\label{patha}
 [\p_o^x(\theta^*_o)]_i=(1-\alpha^x)x_{io}+\alpha^xx^m_i, 
 \end{equation} 
 and hence
 the derivative of $[\p^x({\theta^*})]_i$ with respect to  $x_{io}$ is positive. 
The proof for $[\p^y({\theta^*})]_r$ is analogous.
\end{proof}
\begin{pozn}
It is easy to see that the GS model with the GS range direction does not meet the property \ref{P8} over the trivial technology $\T$. By Lemma \ref{dir} the path generated by directions \eqref{directions} passes through $(\x^{\m}, \y^{\M})$ at $\theta_{\m} $ and since this point belongs to $\T$, it holds $\theta_{\m}=\theta_o^*$. Hence, the efficiency score $\theta_o^*$ of each unit $(\x_o, \y_o)\in\T\setminus\{ (\x^{\m}, \y^{\M})\}$ is equal to the same value $\theta_{\m}$. Note that this is consistent with the results of Section \ref{SS3.3}: since for any $(\x_o, \y_o)\in\T\setminus\{ (\x^{\m}, \y^{\M})\}$ the corresponding projection $\p(\theta^*_o)$ is equal to the same value $\p(\theta^*_o)=\p(\theta_{\m})$,  the path-flow mapping $\p(\theta^*_o)$  is not strictly monotone on $\T$ at $\theta_o^*$.  
\end{pozn}

\section{Numerical examples}\label{S5}

This section documents the ability of individual path-based DEA models to project  units onto the strongly efficient frontier in different real datasets. We consider two basic settings for the function $\psi$: linear, which leads to the DDF-g
models, and hyperbolic, which leads to the HDF-g models. As direction vectors, we have selected vectors from groups analysed by \cite{halicka2024unified}; their description can be found in Table \ref{Tab:choices}.

\begin{table}[H]
\centering
\begin{tabular}{l c c c}
 \hline\\[-1.8ex]
 Notation  & $\g_o^x$  & $\g_o^y$ & Reference\\[0.8ex]
 \hline\\[-1.8ex]
 \newtag{(G1)}{G1} & $\vert\x_o\vert$  & $\vert\y_o\vert$ & \cite{kerstens2011negative} \\
 \newtag{(G2)}{G2} & $\delta^x(\x_o- \x^{\m})$ & $\delta^y (\y^{\M}-\y_o)$ &\cite{halicka2024unified} \\
 \newtag{(G3)}{G3}&$ \x^{\M}- \x^{\m}$&$ \y^{\M}- \y^{\m}$&\cite{portela.al.04}\\
 \newtag{(G4)}{G4}&$\vert\x^{ev}\vert$&$\vert\y^{ev}\vert$&\cite{aparicio2013overall}\\
 [0.8ex]
 \hline
 \end{tabular}
\caption{ Choices of directions $g_o$ used in numerical experiments. Here $x^{ev}_i=\frac 1n\sum_j x_{ij}$ and $y^{ev}_r=\frac 1n\sum_j y_{rj}$. The absolute values in directions (G1) and (G4) accommodate the possibility of 
negative data. The directions (G2)  correspond to the GS range directions defined in \eqref{directions}, i.e., $\delta^x=(1-\psi^x(\theta_{min}))^{-1}$, $\delta^y=(\psi^y(\theta_{min})-1)^{-1}$. }
\label{Tab:choices}
\end{table}

The numerical experiments look at 10 real datasets (see Table~\ref{Tab:data}) that have previously been used in various case studies. The descriptive statistics for each dataset can be found in \ref{D}. It should be noted that Dataset 4 contains negative data, and that the two undesirable outputs in Dataset 10 were treated as inputs. Furthermore, the three datasets (i.e., 2, 8, 9) that did not initially contain input or output values were supplemented by input or output values equal to unity for all units, respectively. The datasets vary in the density of efficient units, with the ratio of strongly efficient units to the total number of units ranging from 22\% for Dataset 10 to 57\% for Dataset 6.

\begin{table}[ht!]
\centering
\footnotesize{
\begin{tabular}{cccccl}
 \hline\\[-1.8ex]
Dataset & $\sharp$ total/eff DMUs & density & $\sharp$ inputs & $\sharp$ outputs & source \\
\hline
1 & 20/9  & 45.0\% & 3  & 2  & \cite{sueyoshi2007computational}  \\
2 & 139/32 & 23.0\% &  4  & 0  & \cite{toloo2014finding}  \\
3 & 31/11  & 35.5\% & 3  &  2 & \cite{juo.al.15.profit}   \\
4 & 30/9  & 30.0\% & 2  & 2  &  \cite{tone2020handling} \\
5 & 13/5  & 38.5\% & 1  & 3  & \cite{talluri2000cone}  \\
6 &  28/16 & 57.1\% & 4  & 2  & \cite{ray2008directional}  \\
7 & 31/15  & 48.4\%  & 2  & 3  &  \cite{xiong2019multi} \\
8 & 46/11  & 23.9\% & 0  & 4  & \cite{foroughi2011note}  \\
9 & 15/4  & 26.7\% & 0  & 6  & \cite{liu2011study}  \\
10 & 92/20  & 21.7\% & 3+2  & 1  & \cite{fare2007environmental} \\
\hline
\end{tabular}}
\caption{Basic description of the datasets used in the experiments.}
\label{Tab:data}
\end{table}

We have applied the DDF-g ( $\psi^x(\theta)=\theta$, $\psi^y(\theta)=2-\theta$, $\theta_{\m}=0$) and HDF-g models ( $\psi^x(\theta)=\theta$, $\psi^y(\theta)=\frac{1}{\theta}$,   $\theta_{min}=0.1$), combined with directions (G1), (G2), (G3),  and (G4) to all 10 datasets listed in Table~\ref{Tab:data}. For numerical computations we 
have utilised a Matlab implementation of the CVX modelling system (see \citealp{cvx}, \citealp{grant.boyd.08}) to solve convex programmes. 
 After the efficiency evaluation  using \eqref{general}, we have applied the standard second phase method (described for path-based models by programme (18) in \citealp{halicka2024unified}) to determine the number of units projected onto the strongly efficient frontier. The results are shown  in Table~\ref{Tab:percentage}.
 
\begin{table}[ht]
\centering
\footnotesize{
\begin{tabular}{cllllllll}
 \hline\\[-2.5ex]
                 & \multicolumn{4}{c}{\textbf{DDM-g}}                                                                                                        & \multicolumn{4}{c}{\textbf{HDM-g}}                                                                                                        \\[0.5ex] \hline
\text{Dataset} & \text{(G1)}                      & \text{(G2)}                      & \text{(G3)}                      & \text{(G4)}                      & \text{(G1)}                      & \text{(G2)}                    & \text{(G3)}                      & \text{(G4)}                      \\ \hline \\[-2.5ex]
1  & \bf 50.0  & \bf 50.0  & \bf 50.0  & \bf 50.0  & \bf 50.0  & \bf 50.0  & \bf 50.0  & \bf 50.0  \\
2  & 22.3   & \bf 53.2  & 27.3   & 22.3   & 23.0   & \bf 53.2  & 27.3   & 23.7   \\
3  & \bf 35.5  & \bf 35.5  & \bf 35.5  & \bf 35.5  & \bf35.5  & \bf 35.5  & \bf 35.5  & \bf 35.5  \\
4  & 40.0  & \bf 70.0  & 66.7  & 43.3   & 43.3  & \bf 76.7  & 66.7  & 53.3   \\
5 & 53.9   & \bf 61.5  & 53.9  & \bf 61.5  & 46.2  & \bf 61.5  & 53.9  & \bf 61.5  \\
6  & \bf 57.1  & \bf 57.1  & \bf 57.1  & \bf 57.1  & \bf 57.1  & \bf 57.1  & \bf 57.1  & \bf 57.1  \\
7                & 61.3                          & \bf 67.7  & 58.1                          & 58.1                          & 61.3                          & \bf 64.7  & 58.1                          & 58.1                          \\
8                & 23.9                          & \bf 80.4  & 52.2                          & 23.9                          & 23.9                          & \bf 80.4  & 52.2                          & 23.9                          \\
9                & \bf 26.7  & \bf 26.7  & \bf 26.7  & \bf 26.7  & \bf 26.7  & \bf 26.7  & \bf 26.7  & \bf 26.7  \\
10               & \bf 22.8  & \bf 22.8  & 21.7                          & 21.7                          & 22.8                          & \bf 23.9  & 21.7                          & 22.8                          \\[0.5ex] \hline
\end{tabular}}
\caption{Percentage of units projected onto the strongly efficient frontier $\partial^S\T$ for each model, with maximal value in bold. It can be seen that the GS range direction (G2) always provides the maximal number of strongly efficient projections. }
\label{Tab:percentage}
\end{table}

We can see that four of the  datasets (i.e., 1, 3, 6, 9) are resistant to direction selection --- the same number of units is projected onto the strongly efficient frontier regardless of our choice of directions. A possible explanation for this phenomenon is 
that the units projected onto $\partial^W\T$ are so close to this part of the boundary that even the choice of the direction leading to the ideal point cannot prevent their projections onto $\partial^W\T$. To a lesser extent, this is also true for datasets 5 and 10. In all other datasets, the GS range direction (G2) outperforms, often significantly, the remaining directions.

In our numerical experiment,  
the GS range direction achieves the highest number of strongly efficient projections in each dataset. 
However, as one can see on the example of Datasets 3, 9, and 10, even the the best direction may project majority of the inefficient units  
onto $\partial^W\T$. Their scores do not capture all the inefficiencies. Therefore, it is important to bear this in mind when interpreting the results of path-based models.

\section{Conclusions}\label{S6}
The paper analyses connections among three desirable properties of DEA models: indication, strict monotonicity, and strong efficiency of projections. For a correct interpretation of the results of a model, it is important to know whether the model meets these characteristics. 
A good understanding of the properties allows one to decide whether the units with a score equal to one were correctly identified as strongly efficient, whether the score of any unit from $\T$ captures all sources of inefficiency, and whether the obtained scores of the units allow a fair comparison of inefficient units with each other. If it is known that one of the properties is not met in the model, it may be necessary to perform further analyses or exercise care when interpreting the results of the model.

\begin{figure}[ht!]
    \centering
\adjustbox{scale=0.9,center}{\begin{tikzcd}[sep=large,
cells={nodes={circle, draw, thick, inner xsep=6pt}, inner ysep=6pt, outer xsep=10pt},
every arrow/.append style = {-stealth},arrows={tail arrow, tail arrow tip width=1.3em, tail arrow tail width=0.4em, tail arrow tip length=0.8em,tail arrow joint width=0.4em}    ]
 &[-6ex] |[ellipse,inner xsep=0pt, inner ysep = 2pt, outer ysep = 10pt, outer xsep = 10pt]| 
\parbox{4cm}{\centering \bf\textsf{ideal technology} \\ $+$ \\ \textsf{GS range direction}} 
\arrow[ld, start anchor = {[yshift = 3ex,xshift = -4ex]}, end anchor = {[yshift = 2ex,xshift = -1ex]}] \arrow[rd,start anchor = {[yshift = 3ex,xshift = 4ex]}, end anchor = {[yshift = 2ex,xshift = 1ex]}]  
&[-6ex]      \\[-1ex]
\text{\bf \textcolor{red}{\textsf{MO}}} \ar[rr, shift left=0.2cm, "{\raisebox{0.5em}{\text{\bf \textsf{path-flow continuity}}}}"] \arrow[rd,start anchor = {[yshift = -2ex,xshift = 0ex]},end anchor = {[yshift = 1ex,xshift = 0ex]}] & & 
\text{\bf \mask{\textcolor{red}{\textsf{PR}}}{\textsf{MO}}} 
\ar[ll, shift left=0.2cm, "{\raisebox{-1em}{\text{\bf \textsf{path-flow monotonicity}}}}"] \arrow[ld,start anchor = {[yshift = -2ex,xshift = 0ex]},end anchor = {[yshift = 1ex,xshift = 0ex]}]
\\[1ex]
|[ellipse,inner xsep=0pt, inner ysep = 0pt]| 
\parbox{3.2cm}{\centering \ \\[-0.5ex] \bf\textsf{no point in} $\partial^W\T$ \\ \textsf{has score of 1}\\[-2ex] \ } \arrow[r, yshift=-0.2cm]& \text{\bf \mask{\textcolor{red}{\textsf{ID}}}{\textsf{MO}}} \arrow[l, yshift=0.2cm] \arrow[r]
&  |[ ellipse,inner xsep=0pt]| \parbox{3cm}{\centering \bf \textsf{not all directions}\\ \textsf{are positive}}        
\end{tikzcd}
}
    \caption{Hierarchy of properties for
    the GS model corresponding to a fixed choice of admissible model parameters.}
    \label{F:AC5}
\end{figure}
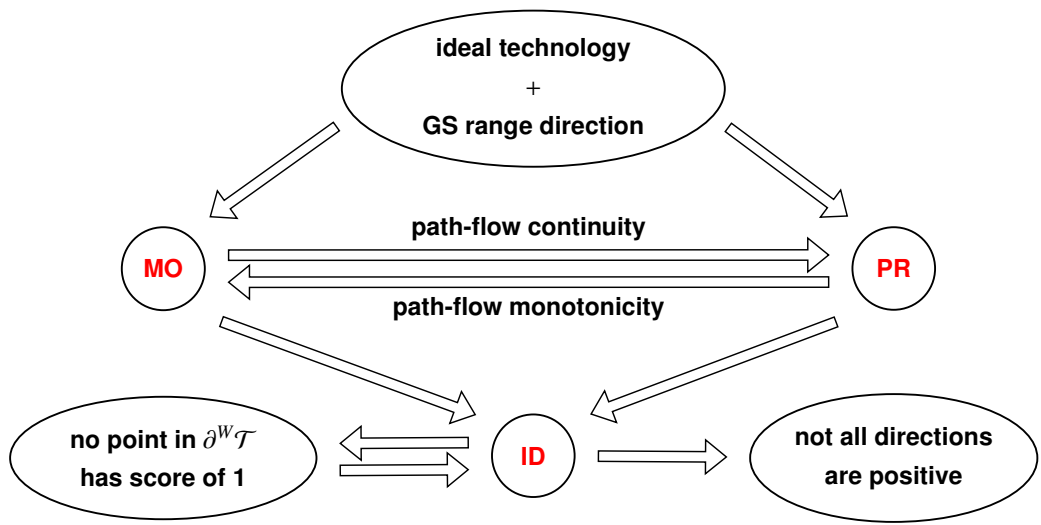

This article focuses on path-based models which are characterised by the fact that individual models, with some exceptions, do not fulfil even one of the three properties. 
Our findings are summarised in Figure~\ref{F:AC5}.

The article shows that the common practice of replacing the property \ref{P3} with the property \ref{P2a} is not justified, in general. This is evidenced by example \ref{ex3.2}, where the model satisfies \ref{P2a} but not \ref{P3}. On the other hand, it was shown (surprisingly to us) that the property \ref{P3} is equivalent to the property \ref{P8} under mild assumptions. 
Nonetheless, the verification of these two properties has a different flavour: while the property \ref{P8} is usually verified analytically for a specific model (i.e., it is proved theoretically as a global property, even outside the class of path-based models), the validation of the property \ref{P3} for a given unit in $\T$ can be performed numerically by applying the standard second phase.

In this paper, we also characterise the cases 
that ensure the three properties are met. Such models are characterised by the GS range directions \eqref{directions} and the data $\X,\Y$ that generate an ideal technology (Definition~\ref{D:idealT}). The latter requirement can be easily verified in practice using Theorem \ref{alternativy}(iv).
Taken as a whole, our analysis yields the following recommendations for DEA practitioners. 
\begin{enumerate}
    \item When using path-based models, it must be assumed that the model does not meet the indication, strict monotonicity, and efficiency of the projections. This means that some units may not be projected onto the strongly efficient boundary. Therefore, if one wishes to know whether the obtained score captures all sources of inefficiency, it is necessary to apply the second-phase to identify the non-included slacks.
    \item To increase the number of units projected onto the strongly efficient frontier, it is advisable to use the GS range directions. These directions ensure that all units are projected onto the strongly efficient frontier in the case of single input, single output data, and more generally for all data configurations that correspond to
    an ideal technology. 
    Although higher-dimensional data typically do not result in an ideal technology, 
    empirical evidence suggests that the GS range directions perform no worse and in some datasets significantly better than their alternatives.
\end{enumerate}

\section*{Acknowledgements }
The authors thank three anonymous reviewers for their valuable comments. 
The research of the first two authors was supported by the APVV-20-0311 project of the Slovak Research and Development Agency and the VEGA 1/0611/21 grant administered jointly by the Scientific Grant Agency of the Ministry of Education, Science, Research and Sport of the Slovak Republic and the Slovak Academy of Sciences.

\appendix
\renewcommand{\thetheorem}{\Alph{section}.\arabic{theorem}}

\section{Polyhedral characteristics of the VRS technology set \texorpdfstring{$\T$}{T}}\label{A}

   As is well known, every polyhedral set can be expressed as the intersection of a limited number of half-spaces (see, e.g., \citealp{klee1959some}, or \citealp{rockafellar.70}).
In the case of the variable returns-to-scale technology set $\T$ defined in \eqref{T},
 the half-spaces must be of the form
\begin{equation}\label{hyperabc}
\HP (\bm{u},\bm{v},\sigma)=\{(\x,\y): \bm{u}^\top\y-\bm{v}^\top\x\le \sigma\}, \quad
\text{where $\bm{u}\in\mathbb{R}^m_{+}$ and $\bm{v}\in\mathbb{R}^s_{+}$.}
\end{equation}
    In order to exclude redundant half-spaces, we  consider only the  \emph{facet-defining}  half-spaces for which 
\begin{equation}\label{facetabc}
\partial \HP (\bm{u},\bm{v},\sigma)\cap\T\equiv \{(\x,\y)\in\T: \bm{u}^\top\y-\bm{v}^\top\x= \sigma\}, 
\end{equation}
is a facet of the polyhedron $\T$.%
\footnote{ The relative boundary $\partial \HP (\bm{u},\bm{v},\sigma)$
of the facet defining half-space $\HP (\bm{u},\bm{v},\sigma)$ is a supporting hyperplane to $\T$ at any point of the corresponding facet.}
 We also introduce notation $\HP_i^x$ and $\HP_i^y$ for
 the special half-spaces 
\begin{equation}\label{eq:HPir}
\begin{split}
\HP_i^x:=
\{(\x,\y): x_i\ge x_i^{\m}\},\qquad
\HP_r^y:=
\{(\x,\y): y_r\le y_r^{\M}\},
\end{split}
\end{equation} 
and denote the corresponding facets by
 \begin{equation}\label{F}
\F_i^x=\{(\x,\y)\in\T: x_i=x_i^{\m}\}, \qquad
\F_r^y=\{(\x,\y)\in\T: y_r=y_r^{\M}\},
\end{equation}
for $i=1,\dots,m$ and $r=1,\dots,s$. 
The following lemma lists some known results about polyhedral sets applied to the facial structure of general VRS technologies defined by \eqref{T}. These results are helpful in the proof of Theorem~\ref{alternativy}. 

\begin{lemma} \label{generalT} Any technology set $\T$ defined by \eqref{T} is the intersection of a finite number of facet-defining half-spaces  
$\HP (\uu,\vv,\sigma)$ with $\uu\ge 0$ and $\vv\ge 0$ and its boundary is the union of the corresponding facets. Moreover, the following statements hold.
\begin{enumerate}[(a)]
\item\label{generalT.a} Facets with $\uu>\0$ and $\vv>\0$ are bounded and their union forms $\partial^S\T$.
\item\label{generalT.b} Facets whose vector $(\uu,\vv)$ contains at least one zero component are unbounded and their union is the closure of $\partial^W\T$. 
\item\label{generalT.c} Each of the special half-spaces $\HP_i^x,\HP_r^y$ in \eqref{eq:HPir} is facet-defining. Furthermore, the corresponding facets $\F^x_i$, $\F^y_r$  are unbounded. 
\item\label{generalT.d} Every point $(\x,\y)\in \T\setminus(\cup_{i=1}^m \F_i^x)\cup (\cup_{r=1}^s \F_r^y)$ satisfies $\x>\x^{\m}$ and $\y<\y^{\M}$. 
\end{enumerate}
\end{lemma}

In line with \cite{davtalab2015characterizing}, we refer to the facets in \ref{generalT.b}, including those in \ref{generalT.c},  as \emph{weak facets}.

\section{Characterisation of ideal technology sets}\label{B}

The next theorem formulates six equivalent characterisations of ideal technology sets, of which the property \ref{ek.iii} was previously used as the definition (see \eqref{PW}  in Subsection \ref{ss42}).
By comparing the properties of the general technology set $\T$ presented in Lemma~\ref{generalT} with the characterisation in Theorem~\ref{alternativy} \ref{ek.i}, we see that ideal technologies are those sets whose weak facets are precisely the special facets $\F_i^x$ and $\F_r^y$.

\begin{theorem}\label{alternativy} For a technology set $\T$ defined by \eqref{T}, the following are equivalent.
\begin{enumerate}[(i)]
\item\label{ek.i}  $\T$
is the intersection of the $m+s$ half-spaces $\HP_i^x$ and $\HP_r^y$
and a finite number of facet defining half-spaces  
$\HP (\bm{u},\bm{v},\sigma)$ with positive
 $\bm{u}$ and $\bm{v}$.
\item\label{ek.ii.bis} The closure of the weakly efficient boundary can be represented as $\cl\partial^W\T=(\cup_{i=1}^m \F_i^x)\cup (\cup_{r=1}^s \F_r^y)$. 
\item\label{ek.ii} 
 The weakly efficient boundary can be represented as $\partial^W\T=(\cup_{i=1}^m \F_i^x)\cup (\cup_{r=1}^s \F_r^y)\setminus\partial^S\T$. 
 \item\label{ek.v} 
For each of the $m+s$ input--output coordinates, there is a generating unit DMU$_j$, $j\in\{1,\ldots,n\}$, which coincides with the ideal point except perhaps in this one coordinate. More formally, for each $i'\in \{1,\ldots,m\}$ there exists a $j\in\{1,\dots,n\}$ such that $x_{ij}=x^{\m}_i$ $\forall$  $ i\neq i'$ and $\y_j=\y^{\M}$; and for each $s'\in \{1,...,r\}$ there exists a $j\in\{1,\dots,n\}$ such that $\x_j=\x^{\m}$ and $y_{sj}=y^{\M}_s$, $\forall$~$s\neq s'$.
\item\label{ek.iv} 
For each of the $m+s$ input--output coordinates, there is a point in the technology set  that coincides with the ideal point  except perhaps in this one coordinate.
More formally, for each $i\in \{1,\ldots,m\}$ there exists a
$\delta^x\ge 0,$ such that $(\x^\m, \y^\M)+\delta^x (\e_i^x,0_s)\in \T$ and for each $r\in \{1,\ldots,s\}$ there exists a
$\delta^y\ge 0,$ such that $(\x^\m, \y^\M)-\delta (0_m, \e_r^y)\in \T$. 
\item\label{ek.iii}  Each point on the weakly efficient frontier coincides with the ideal point in at least one component. More formally, for each  $(\x,\y)\in \partial^W \T$ there exists a $ i\in \{1,...,m\}: \ x_i=x_i^\m \ \hbox{or} \
 r\in \{1,...,s\}: \ y_r=y_r^\M. $
\end{enumerate}
\end{theorem}
\begin{proof}  The proof scheme is shown in Figure \ref{F:proof}.
\begin{figure}[t]
    \centering
\newlength{\edgelentgh}
\setlength{\edgelentgh}{2cm}
\begin{tikzcd}[row sep={0cm,between origins},column sep={0cm,between origins}] 
&[.31\edgelentgh]
&[.4\edgelentgh]
\ref{ek.i}\arrow[d,Leftrightarrow] 
&[.4\edgelentgh] 
&[.31\edgelentgh]
\\[.68\edgelentgh]
&&\ref{ek.ii.bis}\arrow[drr,Rightarrow]&& 
\\[.56\edgelentgh]
\ref{ek.iii}\arrow[urr,Rightarrow] 
&&&& 
\ref{ek.ii} \arrow[dl,Rightarrow] 
\\[.65\edgelentgh] 
& \ref{ek.iv} \arrow[lu,Rightarrow] && 
\ref{ek.v} \arrow[ll,Rightarrow] & 
\end{tikzcd}
    \caption{Scheme of proof for Theorem~\ref{alternativy}.}
    \label{F:proof}
\end{figure}

\begin{itemize}[labelwidth =\widthof{(ii) $\Rightarrow$ (iii)},leftmargin = !]
 \item[\ref{ek.i} $\Rightarrow$ \ref{ek.ii.bis}]  Assume that \ref{ek.i} holds. Then \ref{ek.ii.bis} follows from Lemma~\ref{generalT}. 
\item[\ref{ek.ii.bis} $\Rightarrow$  \ref{ek.i}] Assume that \ref{ek.ii.bis} holds and \ref{ek.i} does not. Lemma~A.1 then yields a facet generating half-space $\hat\HP$ different from those in \eqref{eq:HPir} whose vector $(a,b)$ is non-positive. The corresponding facet therefore lies in the closure of $\partial^W\T$ in contradiction to \ref{ek.ii.bis}.
\item[\ref{ek.ii.bis} $\Rightarrow$ \ref{ek.ii}]  
Since $\cl \partial^W\T\subset\partial\T$, the claim follows.
\item [\ref{ek.ii} $\Rightarrow$  \ref{ek.v}]  The claim \ref{ek.ii} implies that $\T$ contains $m+s$ facets $\F_i, \F_r$. These facets are mutually orthogonal and, hence, the intersections of any $m+s-1$ of these facets form edges (one-dimensional faces) of $\T$. Each of these edges contains a vertex of the form described in \ref{ek.v}. Since the vertices of $\T$ belong to the set of units generating $\T$, the assertion \ref{ek.v} holds.
\item [\ref{ek.v} $\Rightarrow$  \ref{ek.iv}] This implication is trivial.
\item [\ref{ek.iv} $\Rightarrow$ \ref{ek.iii}] 
 Assume by contradiction that \ref{ek.iv} holds and that there exists $(\x_o, \y_o)\in \partial^W\T$ such that 
$(\x^\m, \y^\M)\succ (\x_o, \y_o)$. Since $(\x_o, \y_o)$ is weakly but not strongly efficient, due to free disposability, there exists a unit in $\T$ dominating $(\x_o, \y_o)$ in exactly one component. Without loss of generality, assume that it is the $i$-th input component, i.e., one has  $(\x_o, \y_o)-\gamma (\e_i^x,0_m)\in \T$ for some $\gamma >0$. Note that \ref{ek.iv} implies
$(\x^\m, \y^\M)+\delta^x (\e_i^x,0_s)\in \T$ for some $\delta^x\ge 0$. For $\lambda\in (0,1)$ define
$$
(\x(\lambda), \y(\lambda)):=\lambda [(\x^\m, \y^\M)+\delta (\e_i^x,0_s)]+
(1-\lambda)[(\x_o, \y_o)-\gamma (\e_i^x,0_m)].
$$
Clearly $(\x(\lambda), \y(\lambda))\in\T$ for all $\lambda\in (0,1)$. For $\lambda$ sufficiently small, we have 
$$(\x(\lambda), \y(\lambda))\succ (\x_o, \y_o),$$ which contradicts $(\x_o, \y_o)\in \partial^W\T$.

\item [\ref{ek.iii} $\Rightarrow$  \ref{ek.ii.bis}]
Assume that \ref{ek.iii} holds and \ref{ek.ii.bis} does not. Then there exists a point $(\hat\x, \hat\y)$ in $\partial^W\T$ that does not belong to $(\cup_{i=1}^m \F_i^x)\cup (\cup_{r=1}^s \F_r^y)$. Lemma~\ref{generalT}\ref{generalT.d} yields that $\hat\x>\x^{\m}$ and $\hat\y<\y^{\M}$ in contradiction to \ref{ek.iii}. \qedhere 
\end{itemize}
 \end{proof}
 
\section{Proof of Lemma~\ref{AcorPW}}\label{C}

Trivially, 
$\mathcal \F_{IR}=(\cap_{i\in I}\F_i^x)\cap(\cap_{r\in R}\F_r^y)$ (see definitions in \eqref{F}).   
 Since the normals of the facets $\F_i$ and $\F_r$  are linearly independent, $ \F_{IR}$ is the $(m+s-\vert I\vert-\vert R\vert)$ - dimensional face of $\T$. The (relative) boundary of each face consists of faces of lower dimensions. Theorem~\ref{alternativy}\ref{ek.i} implies that the relative boundary of $\mathcal \F_{IR}$ consists of two types of $(m+s-\vert I\vert-\vert R\vert-1 )$ - dimensional faces: the first one is expressed as the intersections of $\mathcal \F_{IR}$ with one of the facets $F_i$, $i\notin I$ or  $F_r$, $r\notin R$;   the second one is expressed as the intersections of $\mathcal \F_{IR}$ with one of the faces generated by $\HP (\bm{a},\bm{b},\sigma)$ where $\ab>0$, $\bb>0$.

From the definition of $I$ and $R$ it follows that $(\x,\y)$  does not belongs to  any of the faces 
of the first type. Therefore, $(\x,\y)$ is either an interior point of $\mathcal \F_{IR}$ or it belongs to a face of type $\F (\bm{a},\bm{b},\sigma)$, which due to the positiveness of the vectors $\bm{a},\bm{b}$ belongs to the strongly efficient frontier of $\T$. The lemma is proved.

\setcounter{table}{0}
\section{Descriptive statistics of datasets used in Section \ref{S5}}\label{D}

\begin{table}[!ht]
 \centering
 \footnotesize{
 \begin{tabular}{lrrrrr}
  \hline\\[-2.8ex]
          & $x_1$ & $x_2$ & $x_3$ & $y_1$ & $y_2$ \\
          \hline
min       &   27.00 & 115.00 & 2584.00 & 5765.00 &  3286.00\\
max       & 1043.00 & 436.00 & 17412.00 &  223340.00  &  39653.00 \\
mean      & 403.10 & 248.80 & 8170.90 & 81723.00 & 14413.40 \\
std. dev. & 381.20&  96.37 & 5130.93 & 82604.64 & 11747.93\\
\hline
\end{tabular}}
\caption{\textbf{Dataset 1} (\citealp{sueyoshi2007computational}) consists of data of 20 Japanese banks, with three inputs (total amount of capital, number of
offices, and number of employees) and two outputs (total profit and total amount of deposits).}
\end{table}

\begin{table}[ht]
 \centering 
 \footnotesize{
 \begin{tabular}{lrrrr}
  \hline\\[-2.8ex]
          & $x_1$ & $x_2$ & $x_3$ & $x_4$ \\
          \hline
min       & 0.00 &  6521.00 &  0.000 & 1.192 \\  
max       & 508203.00 &  36062.00   & 1213.000   & 1.804 \\
mean      & 225635.85 &  15382.37 &   976.281   & 1.413 \\
std. dev. & 152560.30  & 7652.33 &  480.586  & 0.151 \\
\hline
    \end{tabular}}
\caption{\textbf{Dataset 2} (\citealp{toloo2014finding}) consists of data of 139 different alternatives for long-term asset financing at Czech banks and leasing companies, with four inputs (downpayment, annuities, other fees, and bank loan coefficient) and no outputs.}
\end{table}

\begin{table}[ht]
 \centering 
 \footnotesize{
 \begin{tabular}{lrrrrr}
  \hline\\[-2.8ex]
          & $x_1$ & $x_2$ & $x_3$ & $y_1$ & $y_2$ \\
          \hline
min       &  24930.00  &  204.00 &  489.00 & 2203.00  &     60552.00\\
max       & 2781454.00  & 8804.00 & 78284.00 &   431255.00  &   1912741.00\\
mean      &  710482.84 & 3811.26 &   14003.23 &  134856.35 &  560873.87\\
std. dev. & 684169.90 &  2608.07 &  15703.33 & 135245.53 & 537324.60 \\
\hline
    \end{tabular}}
\caption{\textbf{Dataset 3} (\citealp{juo.al.15.profit}) consists of data of 31 banks operating in Taiwan, with three inputs (financial funds, labor, and physical capital) and two outputs (financial investments and loans).}
\end{table}

\begin{table}[ht]
 \centering
 \footnotesize{
 \begin{tabular}{lrrrr}
  \hline\\[-2.8ex]
          & $x_1$ & $x_2$ & $y_1$ & $y_2$ \\
          \hline
min       &  418242.00  &  1319.0 & -561965.00 & -0.25830 \\
max       & 16655569.00 & 634436.0 & 3092358.00 & 3.53714 \\
mean      & 3705789.57 &  102759.1 & 4.26691.37 &  0.32565\\
std. dev. &  3807894.37 &  135690.5 & 653124.35 &   0.79420 \\
\hline
\end{tabular}}
\caption{\textbf{Dataset 4} (\citealp{tone2020handling}) consists of data of 30 Taiwanese electrical machinery firms with two inputs (cost of sales and R\&D expenses) and two outputs (net income and return).}
\end{table}

\begin{table}[ht]
 \centering 
 \footnotesize{
 \begin{tabular}{lrrrr}
  \hline\\[-2.8ex]
          & $x_1$ & $y_1$ & $y_2$ & $y_3$ \\
          \hline
min       &  4.0  &  0.025 & 6.00 & 0.50 \\
max       & 12.5 & 1.250 & 115.00 & 2.90 \\
mean      & 7.6 & 0.420 & 44.92 & 1.22\\
std. dev. & 2.6 & 0.450 & 35.16  & 0.73 \\
\hline
\end{tabular}}
\caption{\textbf{Dataset 5} (\citealp{talluri2000cone}) consists of data of 13 industrial robots with single input (cost) and three outputs (repeatability, load capacity,
and velocity).}
\end{table}

\begin{table}[H]
 \centering 
 \footnotesize{
 \begin{tabular}{lrrrrrr}
  \hline\\[-2.8ex]
          & $x_1$ & $x_2$ & $x_3$ & $x_4$ & $y_1$ & $y_2$ \\
          \hline
min       &  4067.00 &  62.00 & 241.00 &  587.00 & 2943.00 &  65.00 \\
max       & 80627.00 &  2381.00 & 5678.00 & 18624.00 & 133796.00  &  5346.00 \\
mean      &  28690.43 & 838.18 & 2140.25 &  7075.04 & 4.3555.86 & 1298.61 \\
std. dev. & 20279.43 & 653.57 &  1442.54 &  5352.85 & 3.4877.86 &   1183.87 \\
\hline
\end{tabular}}
\caption{\textbf{Dataset 6} (\citealp{ray2008directional}) consists of data of 28 international airlines (year 1990)
with four inputs (number of employees, fuel, operating and maintenance expenses, capital) and two outputs (passenger-kilometers flown, freight tonne-kilometers flown).}
\end{table}

\begin{table}[H]
 \centering 
 \footnotesize{
 \begin{tabular}{lrrrrr}
  \hline\\[-2.8ex]
          & $x_1$ & $x_2$ &  $y_1$ & $y_2$ & $y_3$ \\
          \hline
min       & 1.0  & 367.70 & 31.26 &  39.78 &  124.63 \\
max       & 65.9 & 3738.00 & 3209.72 &  1887.45 &  21801.65 \\
mean      & 22.7 & 1609.20 &  1132.71 &  0.738.08 & 5383.30 \\
std. dev. & 18.7 & 896.72 & 812.36 &   501.28 &  5357.41 \\
\hline
\end{tabular}}
\caption{\textbf{Dataset 7} (\citealp{liu2011study}) consists of data of 31 Chinese transportation sectors, with two inputs (labour,capital) and three outputs (GDP by transportation sector, passenger turnover volume, freight turnover volume).}
\end{table}

\begin{table}[H]
 \centering 
 \footnotesize{
 \begin{tabular}{lrrrr}
  \hline\\[-2.8ex]
          & $y_1$ & $y_2$ & $y_3$ & $y_4$ \\
          \hline
min       & 1.18 &  10.82  & 163.0 &  11.29  \\
max       & 3.87 & 55.65 & 718.0 &  139.02 \\
mean      & 1.81 &  25.05 & 474.7 &  56.52 \\
std. dev. & 0.72  & 11.21 &  138.99 &  31.12\\
\hline
\end{tabular}}
\caption{\textbf{Dataset 8} (\citealp{foroughi2011note}) consists of data of 46 association rules, with no inputs  and four outputs (support, confidence, itemset value, cross-selling profit).}
\end{table}

\begin{table}[H]
 \centering 
 \footnotesize{
 \begin{tabular}{lrrrrrr}
  \hline\\[-2.8ex]
          & $y_1$ & $y_2$ & $y_3$ & $y_4$ & $y_5$ & $y_6$ \\
          \hline
min       & 6.88 &  12.06 &   1.15 &   2.68 &  11.22 &  15.97 \\
max       & 100.00 &  100.00 &   100.00 &  100.00 &  100.00 &   100.00 \\
mean      & 36.56 &  47.80 &  16.10 &  19.08 &  43.67 &  48.54 \\
std. dev. & 32.09 & 30.15 &  26.25 &  25.51 & 22.38 &  26.82 \\
\hline
\end{tabular}}
\caption{\textbf{Dataset 9} (\citealp{liu2011study}) consists of data of 15 research institutions, with no inputs  and six outputs (SCI   pub. per research staff, SCI pub. per research expenditures, high quality papers per research staff, high quality papers per research expenditures, external research funding, graduation students enrolled).}
\end{table}

\begin{table}[H]
 \centering 
 \footnotesize{
 \begin{tabular}{lrrrrrr}
  \hline\\[-2.8ex]
          & $x_1$ & $x_2$ & $x_3$ & $b_1$ & $b_2$ & $y_1$ \\
          \hline
min       &  39349184.0 &  39.0 &  1892407000000  & 1293.20  & 423.05 & 166616000 \\
max       & 750024804.0 &  535.0 &  175688406305600.0 &  252344.60 &   72524.10 & 18212069000 \\
mean      &  240000014.7 &  185.2 &   47104363123213.0 &  40745.19 &  17494.02 & 4686524843 \\
std. dev. & 146352514.9 &  110.9  & 39982238625707.7 &  48244.78  & 16190\\
\hline
\end{tabular}}
\caption{\textbf{Dataset 10} (\citealp{fare2007environmental}) consists of data of 92 power plants, with three inputs (capital, labour, heat), two undesirable outputs ($SO_2$ and $NO_x$ emissions), which are handled as inputs,  and single output (energy).}
\end{table}


\begin{thebibliography}{37}
\expandafter\ifx\csname natexlab\endcsname\relax\def\natexlab#1{#1}\fi
\providecommand{\url}[1]{\texttt{#1}}
\providecommand{\href}[2]{#2}
\providecommand{\path}[1]{#1}
\providecommand{\DOIprefix}{doi:}
\providecommand{\ArXivprefix}{arXiv:}
\providecommand{\URLprefix}{URL: }
\providecommand{\Pubmedprefix}{pmid:}
\providecommand{\doi}[1]{\href{http://dx.doi.org/#1}{\path{#1}}}
\providecommand{\Pubmed}[1]{\href{pmid:#1}{\path{#1}}}
\providecommand{\bibinfo}[2]{#2}
\ifx\xfnm\relax \def\xfnm[#1]{\unskip,\space#1}\fi
\bibitem[{Aparicio et~al.(2013)Aparicio, Pastor and Ray}]{aparicio2013overall}
\bibinfo{author}{Aparicio, J.}, \bibinfo{author}{Pastor, J.T.},
  \bibinfo{author}{Ray, S.C.}, \bibinfo{year}{2013}.
\newblock \bibinfo{title}{An overall measure of technical inefficiency at the
  firm and at the industry level: {T}he ‘lost profit on outlay’}.
\newblock \bibinfo{journal}{European Journal of Operational Research}
  \bibinfo{volume}{226}, \bibinfo{pages}{154--162}.
\bibitem[{Banker et~al.(1984)Banker, Charnes and Cooper}]{banker.al.84}
\bibinfo{author}{Banker, R.D.}, \bibinfo{author}{Charnes, A.},
  \bibinfo{author}{Cooper, W.W.}, \bibinfo{year}{1984}.
\newblock \bibinfo{title}{Some models for estimating technical and scale
  inefficiencies in data envelopment analysis}.
\newblock \bibinfo{journal}{Management Science} \bibinfo{volume}{30},
  \bibinfo{pages}{1078--1092}.
\newblock \DOIprefix\doi{10.1287/mnsc.30.9.1078}.
\bibitem[{Chambers et~al.(1996)Chambers, Chung and
  F{\"a}re}]{chambers1996benefit}
\bibinfo{author}{Chambers, R.G.}, \bibinfo{author}{Chung, Y.},
  \bibinfo{author}{F{\"a}re, R.}, \bibinfo{year}{1996}.
\newblock \bibinfo{title}{Benefit and distance functions}.
\newblock \bibinfo{journal}{Journal of Economic Theory} \bibinfo{volume}{70},
  \bibinfo{pages}{407--419}.
\bibitem[{Chambers et~al.(1998)Chambers, Chung and F{\"a}re}]{chambers.al.98}
\bibinfo{author}{Chambers, R.G.}, \bibinfo{author}{Chung, Y.},
  \bibinfo{author}{F{\"a}re, R.}, \bibinfo{year}{1998}.
\newblock \bibinfo{title}{Profit, directional distance functions, and
  {N}erlovian efficiency}.
\newblock \bibinfo{journal}{Journal of Optimization Theory and Applications}
  \bibinfo{volume}{98}, \bibinfo{pages}{351--364}.
\newblock \DOIprefix\doi{10.1023/A:1022637501082}.
\bibitem[{Charnes et~al.(1985)Charnes, Cooper, Golany, Seiford and
  Stutz}]{charnes.al.85}
\bibinfo{author}{Charnes, A.}, \bibinfo{author}{Cooper, W.W.},
  \bibinfo{author}{Golany, B.}, \bibinfo{author}{Seiford, L.},
  \bibinfo{author}{Stutz, J.}, \bibinfo{year}{1985}.
\newblock \bibinfo{title}{Foundations of data envelopment analysis for
  {P}areto-{K}oopmans efficient empirical production functions}.
\newblock \bibinfo{journal}{Journal of Econometrics} \bibinfo{volume}{30},
  \bibinfo{pages}{91--107}.
\newblock \DOIprefix\doi{10.1016/0304-4076(85)90133-2}.
\bibitem[{Chavas and Cox(1999)}]{chavas1999generalized}
\bibinfo{author}{Chavas, J.P.}, \bibinfo{author}{Cox, T.L.},
  \bibinfo{year}{1999}.
\newblock \bibinfo{title}{A generalized distance function and the analysis of
  production efficiency}.
\newblock \bibinfo{journal}{Southern Economic Journal} \bibinfo{volume}{66},
  \bibinfo{pages}{294--318}.
\bibitem[{Cooper et~al.(1999)Cooper, Park and Pastor}]{cooper.al.99}
\bibinfo{author}{Cooper, W.W.}, \bibinfo{author}{Park, K.S.},
  \bibinfo{author}{Pastor, J.T.}, \bibinfo{year}{1999}.
\newblock \bibinfo{title}{{RAM}: A range adjusted measure of inefficiency for
  use with additive models, and relations to other models and measures in
  {DEA}}.
\newblock \bibinfo{journal}{Journal of Productivity Analysis}
  \bibinfo{volume}{11}, \bibinfo{pages}{5--42}.
\newblock \DOIprefix\doi{10.1023/a:1007701304281}.
\bibitem[{Cooper et~al.(2011)Cooper, Pastor, Borras, Aparicio and
  Pastor}]{cooper2011bam}
\bibinfo{author}{Cooper, W.W.}, \bibinfo{author}{Pastor, J.T.},
  \bibinfo{author}{Borras, F.}, \bibinfo{author}{Aparicio, J.},
  \bibinfo{author}{Pastor, D.}, \bibinfo{year}{2011}.
\newblock \bibinfo{title}{{BAM}: A bounded adjusted measure of efficiency for
  use with bounded additive models}.
\newblock \bibinfo{journal}{Journal of Productivity Analysis}
  \bibinfo{volume}{35}, \bibinfo{pages}{85--94}.
\bibitem[{Cooper et~al.(2007)Cooper, Seiford and Tone}]{cooper.al.07}
\bibinfo{author}{Cooper, W.W.}, \bibinfo{author}{Seiford, L.M.},
  \bibinfo{author}{Tone, K.}, \bibinfo{year}{2007}.
\newblock \bibinfo{title}{Data Envelopment Analysis: {A} Comprehensive Text
  with Models, Applications, References and DEA-Solver Software}.
\newblock \bibinfo{publisher}{Springer {US}}.
\newblock \DOIprefix\doi{10.1007/978-0-387-45283-8}.
\bibitem[{Davtalab-Olyaie et~al.(2015)Davtalab-Olyaie, Roshdi, Partovi~Nia and
  Asgharian}]{davtalab2015characterizing}
\bibinfo{author}{Davtalab-Olyaie, M.}, \bibinfo{author}{Roshdi, I.},
  \bibinfo{author}{Partovi~Nia, V.}, \bibinfo{author}{Asgharian, M.},
  \bibinfo{year}{2015}.
\newblock \bibinfo{title}{On characterizing full dimensional weak facets in
  {DEA} with variable returns to scale technology}.
\newblock \bibinfo{journal}{Optimization} \bibinfo{volume}{64},
  \bibinfo{pages}{2455--2476}.
\newblock \DOIprefix\doi{10.1080/02331934.2014.917305}.
\bibitem[{F\"{a}re et~al.(1985)F\"{a}re, Grosskopf and Lovell}]{fare.al.85}
\bibinfo{author}{F\"{a}re, R.}, \bibinfo{author}{Grosskopf, S.},
  \bibinfo{author}{Lovell, C.A.K.}, \bibinfo{year}{1985}.
\newblock \bibinfo{title}{The Measurement of Efficiency of Production}.
\newblock \bibinfo{publisher}{Springer Netherlands}.
\newblock \DOIprefix\doi{10.1007/978-94-015-7721-2}.
\bibitem[{F{\"a}re et~al.(2007)F{\"a}re, Grosskopf and
  Pasurka~Jr}]{fare2007environmental}
\bibinfo{author}{F{\"a}re, R.}, \bibinfo{author}{Grosskopf, S.},
  \bibinfo{author}{Pasurka~Jr, C.A.}, \bibinfo{year}{2007}.
\newblock \bibinfo{title}{Environmental production functions and environmental
  directional distance functions}.
\newblock \bibinfo{journal}{Energy} \bibinfo{volume}{32},
  \bibinfo{pages}{1055--1066}.
\bibitem[{F\"{a}re and Lovell(1978)}]{fare.lovell.78}
\bibinfo{author}{F\"{a}re, R.}, \bibinfo{author}{Lovell, C.A.K.},
  \bibinfo{year}{1978}.
\newblock \bibinfo{title}{Measuring the technical efficiency of production}.
\newblock \bibinfo{journal}{Journal of Economic Theory} \bibinfo{volume}{19},
  \bibinfo{pages}{150--162}.
\newblock \DOIprefix\doi{10.1016/0022-0531(78)90060-1}.
\bibitem[{Foroughi(2011)}]{foroughi2011note}
\bibinfo{author}{Foroughi, A.}, \bibinfo{year}{2011}.
\newblock \bibinfo{title}{A note on “a new method for ranking discovered
  rules from data mining by dea”, and a full ranking approach}.
\newblock \bibinfo{journal}{Expert Systems with Applications}
  \bibinfo{volume}{38}, \bibinfo{pages}{12913--12916}.
\bibitem[{Grant and Boyd(2008)}]{grant.boyd.08}
\bibinfo{author}{Grant, M.}, \bibinfo{author}{Boyd, S.}, \bibinfo{year}{2008}.
\newblock \bibinfo{title}{Graph implementations for nonsmooth convex programs},
  in: \bibinfo{editor}{Blondel, V.}, \bibinfo{editor}{Boyd, S.},
  \bibinfo{editor}{Kimura, H.} (Eds.), \bibinfo{booktitle}{Recent Advances in
  Learning and Control}. \bibinfo{publisher}{Springer-Verlag Limited}. Lecture
  Notes in Control and Information Sciences, pp. \bibinfo{pages}{95--110}.
\newblock \bibinfo{note}{\url{http://stanford.edu/~boyd/graph_dcp.html}}.
\bibitem[{Grant and Boyd(2014)}]{cvx}
\bibinfo{author}{Grant, M.}, \bibinfo{author}{Boyd, S.}, \bibinfo{year}{2014}.
\newblock \bibinfo{title}{{CVX}: Matlab software for disciplined convex
  programming, version 2.1}.
\newblock \bibinfo{howpublished}{\url{http://cvxr.com/cvx}}.
\bibitem[{Halick\'{a} and Trnovsk\'{a}(2021)}]{halicka.trnovska.21}
\bibinfo{author}{Halick\'{a}, M.}, \bibinfo{author}{Trnovsk\'{a}, M.},
  \bibinfo{year}{2021}.
\newblock \bibinfo{title}{A unified approach to non-radial graph models in data
  envelopment analysis: {C}ommon features, geometry, and duality}.
\newblock \bibinfo{journal}{European Journal of Operational Research}
  \bibinfo{volume}{289}, \bibinfo{pages}{611--627}.
\newblock \DOIprefix\doi{10.1016/j.ejor.2020.07.019}.
\bibitem[{Halick{\'a} et~al.(2024)Halick{\'a}, Trnovsk{\'a} and
  {\v{C}}ern{\'y}}]{halicka2024unified}
\bibinfo{author}{Halick{\'a}, M.}, \bibinfo{author}{Trnovsk{\'a}, M.},
  \bibinfo{author}{{\v{C}}ern{\'y}, A.}, \bibinfo{year}{2024}.
\newblock \bibinfo{title}{A unified approach to radial, hyperbolic, and
  directional efficiency measurement in data envelopment analysis}.
\newblock \bibinfo{journal}{European Journal of Operational Research}
  \bibinfo{volume}{312}, \bibinfo{pages}{298--314}.
\bibitem[{Juo et~al.(2015)Juo, Fu, Yu and Lin}]{juo.al.15.profit}
\bibinfo{author}{Juo, J.C.}, \bibinfo{author}{Fu, T.T.}, \bibinfo{author}{Yu,
  M.M.}, \bibinfo{author}{Lin, Y.H.}, \bibinfo{year}{2015}.
\newblock \bibinfo{title}{Profit-oriented productivity change}.
\newblock \bibinfo{journal}{Omega} \bibinfo{volume}{57},
  \bibinfo{pages}{176--187}.
\bibitem[{Kerstens and Van~de Woestyne(2011)}]{kerstens2011negative}
\bibinfo{author}{Kerstens, K.}, \bibinfo{author}{Van~de Woestyne, I.},
  \bibinfo{year}{2011}.
\newblock \bibinfo{title}{Negative data in dea: a simple proportional distance
  function approach}.
\newblock \bibinfo{journal}{Journal of the Operational Research Society}
  \bibinfo{volume}{62}, \bibinfo{pages}{1413--1419}.
\bibitem[{Klee(1959)}]{klee1959some}
\bibinfo{author}{Klee, V.}, \bibinfo{year}{1959}.
\newblock \bibinfo{title}{Some characterizations of convex polyhedra}.
\newblock \bibinfo{journal}{Acta Mathematica} \bibinfo{volume}{102},
  \bibinfo{pages}{79--107}.
\newblock \DOIprefix\doi{10.1007/BF02559569}.
\bibitem[{Liu et~al.(2011)Liu, Zhang, Meng, Li and Xu}]{liu2011study}
\bibinfo{author}{Liu, W.}, \bibinfo{author}{Zhang, D.}, \bibinfo{author}{Meng,
  W.}, \bibinfo{author}{Li, X.}, \bibinfo{author}{Xu, F.},
  \bibinfo{year}{2011}.
\newblock \bibinfo{title}{A study of dea models without explicit inputs}.
\newblock \bibinfo{journal}{Omega} \bibinfo{volume}{39},
  \bibinfo{pages}{472--480}.
\bibitem[{Pastor et~al.(2022)Pastor, Aparicio and
  Zof{\'\i}o}]{pastor2022benchmarking}
\bibinfo{author}{Pastor, J.T.}, \bibinfo{author}{Aparicio, J.},
  \bibinfo{author}{Zof{\'\i}o, J.L.}, \bibinfo{year}{2022}.
\newblock \bibinfo{title}{Benchmarking Economic Efficiency: Technical and
  Allocative Fundamentals}.
\newblock International Series in Operations Research and Management Science,
  \bibinfo{publisher}{Springer Cham}.
\bibitem[{Pastor et~al.(1999)Pastor, Ruiz and Sirvent}]{pastor.al.99.enhanced}
\bibinfo{author}{Pastor, J.T.}, \bibinfo{author}{Ruiz, J.L.},
  \bibinfo{author}{Sirvent, I.}, \bibinfo{year}{1999}.
\newblock \bibinfo{title}{An enhanced {DEA} {R}ussell graph efficiency
  measure}.
\newblock \bibinfo{journal}{European Journal of Operational Research}
  \bibinfo{volume}{115}, \bibinfo{pages}{596--607}.
\bibitem[{Portela et~al.(2004)Portela, Thanassoulis and
  Simpson}]{portela.al.04}
\bibinfo{author}{Portela, M.C.A.S.}, \bibinfo{author}{Thanassoulis, E.},
  \bibinfo{author}{Simpson, G.}, \bibinfo{year}{2004}.
\newblock \bibinfo{title}{Negative data in {DEA}: {A} directional distance
  approach applied to bank branches}.
\newblock \bibinfo{journal}{Journal of the Operational Research Society}
  \bibinfo{volume}{55}, \bibinfo{pages}{1111--1121}.
\newblock \DOIprefix\doi{10.1057/palgrave.jors.2601768}.
\bibitem[{Ray(2008)}]{ray2008directional}
\bibinfo{author}{Ray, S.C.}, \bibinfo{year}{2008}.
\newblock \bibinfo{title}{The directional distance function and measurement of
  super-efficiency: an application to airlines data}.
\newblock \bibinfo{journal}{Journal of the Operational Research Society}
  \bibinfo{volume}{59}, \bibinfo{pages}{788--797}.
\bibitem[{Rockafellar(1970)}]{rockafellar.70}
\bibinfo{author}{Rockafellar, R.T.}, \bibinfo{year}{1970}.
\newblock \bibinfo{title}{Convex analysis}.
\newblock Princeton Mathematical Series, No. 28, \bibinfo{publisher}{Princeton
  University Press}, \bibinfo{address}{Princeton, N.J.}
\bibitem[{Russell(1985)}]{russell1985measures}
\bibinfo{author}{Russell, R.R.}, \bibinfo{year}{1985}.
\newblock \bibinfo{title}{Measures of technical efficiency}.
\newblock \bibinfo{journal}{Journal of Economic theory} \bibinfo{volume}{35},
  \bibinfo{pages}{109--126}.
\bibitem[{Russell and Schworm(2011)}]{russell.schworm.11}
\bibinfo{author}{Russell, R.R.}, \bibinfo{author}{Schworm, W.},
  \bibinfo{year}{2011}.
\newblock \bibinfo{title}{Properties of inefficiency indexes on $\langle$input,
  output$\rangle$ space}.
\newblock \bibinfo{journal}{Journal of Productivity Analysis}
  \bibinfo{volume}{36}, \bibinfo{pages}{143--156}.
\newblock \DOIprefix\doi{10.1007/s11123-011-0209-3}.
\bibitem[{Russell and Schworm(2018)}]{russell.schworm.18}
\bibinfo{author}{Russell, R.R.}, \bibinfo{author}{Schworm, W.},
  \bibinfo{year}{2018}.
\newblock \bibinfo{title}{Technological inefficiency indexes: {A} binary
  taxonomy and a generic theorem}.
\newblock \bibinfo{journal}{Journal of Productivity Analysis}
  \bibinfo{volume}{49}, \bibinfo{pages}{17--23}.
\newblock \DOIprefix\doi{10.1007/s11123-017-0518-2}.
\bibitem[{Sueyoshi and Sekitani(2007)}]{sueyoshi2007computational}
\bibinfo{author}{Sueyoshi, T.}, \bibinfo{author}{Sekitani, K.},
  \bibinfo{year}{2007}.
\newblock \bibinfo{title}{Computational strategy for russell measure in dea:
  Second-order cone programming}.
\newblock \bibinfo{journal}{European Journal of Operational Research}
  \bibinfo{volume}{180}, \bibinfo{pages}{459--471}.
\bibitem[{Sueyoshi and Sekitani(2009)}]{sueyoshi.sekitani.09.ejor.196}
\bibinfo{author}{Sueyoshi, T.}, \bibinfo{author}{Sekitani, K.},
  \bibinfo{year}{2009}.
\newblock \bibinfo{title}{An occurrence of multiple projections in {DEA}-based
  measurement of technical efficiency: {T}heoretical comparison among {DEA}
  models from desirable properties}.
\newblock \bibinfo{journal}{European Journal of Operational Research}
  \bibinfo{volume}{196}, \bibinfo{pages}{764--794}.
\newblock \DOIprefix\doi{10.1016/j.ejor.2008.01.045}.
\bibitem[{Talluri and Yoon(2000)}]{talluri2000cone}
\bibinfo{author}{Talluri, S.}, \bibinfo{author}{Yoon, K.P.},
  \bibinfo{year}{2000}.
\newblock \bibinfo{title}{A cone-ratio dea approach for amt justification}.
\newblock \bibinfo{journal}{International Journal of Production Economics}
  \bibinfo{volume}{66}, \bibinfo{pages}{119--129}.
\bibitem[{Toloo and Kresta(2014)}]{toloo2014finding}
\bibinfo{author}{Toloo, M.}, \bibinfo{author}{Kresta, A.},
  \bibinfo{year}{2014}.
\newblock \bibinfo{title}{Finding the best asset financing alternative: A
  dea--weo approach}.
\newblock \bibinfo{journal}{Measurement} \bibinfo{volume}{55},
  \bibinfo{pages}{288--294}.
\bibitem[{Tone(2001)}]{tone2001slacks}
\bibinfo{author}{Tone, K.}, \bibinfo{year}{2001}.
\newblock \bibinfo{title}{A slacks-based measure of efficiency in data
  envelopment analysis}.
\newblock \bibinfo{journal}{European Journal of Operational Research}
  \bibinfo{volume}{130}, \bibinfo{pages}{498--509}.
\bibitem[{Tone et~al.(2020)Tone, Chang and Wu}]{tone2020handling}
\bibinfo{author}{Tone, K.}, \bibinfo{author}{Chang, T.S.}, \bibinfo{author}{Wu,
  C.H.}, \bibinfo{year}{2020}.
\newblock \bibinfo{title}{Handling negative data in slacks-based measure data
  envelopment analysis models}.
\newblock \bibinfo{journal}{European Journal of Operational Research}
  \bibinfo{volume}{282}, \bibinfo{pages}{926--935}.
\bibitem[{Xiong et~al.(2019)Xiong, Chen, An and Wu}]{xiong2019multi}
\bibinfo{author}{Xiong, B.}, \bibinfo{author}{Chen, H.}, \bibinfo{author}{An,
  Q.}, \bibinfo{author}{Wu, J.}, \bibinfo{year}{2019}.
\newblock \bibinfo{title}{A multi-objective distance friction minimization
  model for performance assessment through data envelopment analysis}.
\newblock \bibinfo{journal}{European Journal of Operational Research}
  \bibinfo{volume}{279}, \bibinfo{pages}{132--142}.

\end{thebibliography}
\end{document}